\documentclass[a4paper,12pt,leqno]{amsart}
\usepackage{amsmath}
\usepackage{amsmath}
\usepackage{amssymb}
\usepackage{amsthm}
\usepackage{amsopn}
\usepackage{amscd}
\usepackage{textcomp}
\usepackage[all]{xy}
\usepackage{lscape}
\usepackage{verbatim}

\newtheorem{theo}{Theorem}[section]
\newtheorem{prop}[theo]{Proposition}
\newtheorem{lemme}[theo]{Lemma}

\numberwithin{equation}{section}

\def\aa{{\mathbb A}}
\def\bc{\backslash}

\def\lra{\leftrightarrow}

\numberwithin{equation}{section}
\def\f{{\mathbb F}}

\def\r{{\mathbb R}}

\def\t{{\times}}
\def\z{\mathbb Z}

\def\cc{{\mathbb C}}

\def\lra{\leftrightarrow}

\def\lj{{\bf LJ}}

\def\bc{\backslash}

\def\tr{{{\rm tr}}}

\begin{document}
\title[Global Jacquet-Langlands]
{Global Jacquet-Langlands correspondence for division algebras in characteristic $p$}

\author{A.I.Badulescu}
\address{A.I.Badulescu, Universit\'e Montpellier 2, I3M}
\email{ibadules@math.univ-montp2.fr}
\author{Ph.Roche}
\address{Ph.Roche, Universit\'e Montpellier 2, CNRS, I3M, L2C}
\email{philippe.roche@univ-montp2.fr}
\maketitle

{\bf Abstract:} We prove a full global Jacquet-Langlands correspondence between $GL(n)$ and division algebras over global fields of non zero characteristic. If $D$ is  a central division algebra of dimension $n^2$ over a global field $F$ of non zero characteristic, we prove that there exists an injective map from the set of automorphic representations of $D^\t$ to the set of automorphic square integrable representations of $GL_n(F)$, compatible at all places with the local Jacquet-Langlands correspondence for unitary representations. We characterize the image of the map. As a consequence we get multiplicity one and strong multiplicity one theorems for $D^\t$.

\tableofcontents

\section{Introduction}

We prove the global Jacquet-Langlands correspondence between $GL_n$ over a global field $F$ of characteristic $p$
 and 
$D^\t$ where $D$ is a central division algebra of dimension $n^2$ over $F$. A corollary is 
the multiplicity one and strong multiplicity one theorem for $D^\t$. We then answer two questions asked by Laumon, Rapoport and Stuhler.
The first case of full global Jacquet-Langlands 
correspondence was proved by Jacquet and Langlands \cite{JL}, for $n=2$. This is a monumental work which served as an example for all the other proofs so far. 
 By "full" we mean that there is no condition on the representations to transfer. Notice that "partial" correspondences, say, for automorphic representations which are cuspidal at two places or so, are very useful, but never imply as a corollary the multiplicity one theorem for inner forms.

For $n=3$ and $F$ of zero characteristic a full global correspondence was proved by Flath in \cite{Fl2}.
Then in zero characteristic for $D$ satisfying the additional condition that $D$ is a division algebra at every place where it does not split by Vign\'eras 
(\cite{Vi}, never published) and later by Harris and Taylor (\cite{HT}, Chap. VI). The 
correspondence for every $n$ and 
without condition on $D$ is proved in zero 
characteristic in \cite{Ba} and \cite{BR} (the first paper assumes that $D$ splits at all infinite places, and in 
the second this 
condition is dropped). Only some partial cases of the Jacquet-Langlands correspondences were proved in non zero characteristic, mainly for
 practical purposes 
(need to construct a representation doing this or that), for instance in \cite{He2}, Appendix 2, and \cite{LRS}.

As far as we know, our result  is the first case of full correspondence in non zero characteristic since \cite{JL}. The main ingredients 
not available in the past were the local transfer of all unitary representations and a trace formula in non zero
 characteristic. Laumon and Lafforgue developed
 the trace formula in \cite{Lau} and \cite{Laf}. The formula is not invariant like the one in \cite{AC} so it is more difficult to use. 
 This explains why we had to confine ourselves here only to the case when the inner form is a division algebra.

In the second section we recall the local tools we will use. We are very careful to give reference or full 
arguments for 
results which are "well known" in zero characteristic, but less well known in non zero characteristic. For 
instance we work only 
with functions with support in the regular set (which excludes for example elements whose characteristic
 polynomial is irreducible 
but not separable). The submersion theorem of Harish-Chandra allows one to easily transfer these functions in any 
characteristic. 

In the third section we define the automorphic representations we want to transfer (the discrete series). We use the positive
 characteristic setting 
(as in \cite{Lau} and \cite{Laf}) which is slightly different from the one in zero characteristic but we 
explain how to switch from one 
to another. Then we give the precise claim of our main result.

The fourth section is devoted to the proof. The main ingredient is clearly here the trace formula of Lafforgue. 
Without this non trivial 
result nothing would be possible. We show that the geometric side and the spectral side of the trace 
formula for $GL_n$ take simple form 
when applied to functions coming from $D^\t$.

In the fifth section we give (positive) answer to some questions asked by Laumon, Rapoport and Stuhler in \cite{LRS}. 

The correspondence proved here completes also the proofs of Lubotzky, Samuels and Vishne in \cite{LSV} (see their Remark 1.6).

The first named author was partially supported by the ANR 08-BLAN-0259.

\section{Local}

\subsection{Basic facts}

Let $F$ be a local field and fix an algebraic closure $\bar{F}$ of $F$. Let $D$ be a central division algebra of dimension $d^2$ over $F$.  Let $O_D$ be 
the ring of integers of $D$.  
  
For $r$ any positive integer, we denote $GL_r(D)$ the group of invertible elements of $M_r(D).$ Let $B$ be the   subgroup of upper triangular matrices and 
let standard parabolic subgroups be the parabolic subgroups containing $B$. 
Let $\Delta=\{1,\cdots,r-1\}$, (for $r=1,\Delta=\emptyset$), to any subset $I\subset \Delta$ one associates an ordered partition  $r_I=(r_1,\cdots,r_k)$ of
 $r$  defined by the condition
 $\Delta\setminus I=\{r_1,r_1+r_2,\cdots, r_1+r_2+...+r_{k-1}\}.$ This map is a bijection between the set of subsets of $\Delta$ and the set of ordered
  partitions of $r.$
To any   $I\subset \Delta$ one associates the subgroup
$M_I(D)$ of $GL_r(D)$  which is 
 the group of block diagonal invertible matrices with blocks of size $r_1,r_2,...,r_k$ (the components of the partition $r_I$)  with  coefficients in 
$D,$ the unipotent sub-group 
 $N_I(D)$   which is the group of corresponding upper block triangular matrices with unit matrices on the diagonal  and the associated
  parabolic subgroup $P_I(D)=M_I(D)N_I(D).$  The groups $M_I(D)$ will be called {\bf standard Levi subgroups of $GL_r(D)$}.

If $P$ is associated to $r_I$ consisting of a k-tuple we denote  $\vert P\vert:=k.$
 There is a bijection between the set $\mathcal{P}_0^s$ of standard parabolic subgroups of $GL_r(D)$, the set of ordered partition $R$ 
 of $r,$ and the subsets of $\Delta.$ 
If $P=P_I$ is a standard parabolic subgroup of $GL_r(D)$ we will denote $M_P:=M_I$ the standard Levi component of $P$ and $N_P:=N_I$ its
 unipotent radical.
Two important parabolic subgroups are the one corresponding to $I=\emptyset$ and to $I=\Delta$.
We have $r_{\emptyset}=(1,1,...,1), P_{\emptyset}$ is the standard minimal  parabolic subgroup of  $GL_r(D)$  and $M_0=M_{\emptyset}$ is the
 group $diag(D^{\times},...,D^{\times})$. We have $r_{\Delta}=(r), P_{\Delta}=M_{\Delta}=GL_r(D).$
 
Let $K$ be the maximal compact subgroup $GL_r(O_D)$ of $GL_r(D)$. We endow $GL_r(D)$ with the Haar measure $dg$ such that the volume
 of $K$ is one, and the center $Z$ of $GL_r(D)$ with the Haar measure such that the volume of $Z\cap K$ is one.

\def\O{{\mathcal O}}

\def\tg{\tilde{G}}
 
Set now $G:=GL_r(D)$, $A:=M_r(D)$ and $n:=rd$. The theory of central simple algebras allows one to define the characteristic polynomial $P_g$
for elements $g\in A$ in spite $D$ is non commutative. $P_g$ is a monic polynomial of degree $n$ with coefficients in $F$. It is the main tool for transferring
conjugacy classes between $GL_n(F)$ and its inner forms like $G$.

There are (at least) two ways of defining the characteristic polynomial $P_g$ as we recall hereafter.  For details and proofs see, for example, 
\cite{Pie} chap. 16 and 17. 
It is known by class field theory that the division algebra $D$ contains an unramified extension $E$ of $F$ of degree $d$, with (cyclic) Galois 
group say $Gal(E/F)$, and that
 $A\otimes_F E=M_{n}(E)$ (Corollary 13.3 and Proposition 17.10 \cite{Pie}). This gives an embedding of $A$ into $M_{n}(E)$. If $g$ is an element
  of $A$, the characteristic 
 polynomial $P_g$ of the image of $g$ in $M_{n}(E)$ does not depend on the embedding (by Skolem-Noether theorem). Also, $P_g$ turns out to be 
 stable by all the elements of 
 $Gal(E/F)$, hence $P_g\in F[X]$ and this is the first definition of the characteristic polynomial. An embedding of $A$ in 
$M_n(E)$ preserves the minimal polynomial, so we have that the minimal polynomial of $g$ divides the characteristic polynomial and the roots 
of the characteristic polynomial in 
$\bar{F}$ are also roots of the minimal polynomial.

The other way of defining $P_g$ is the following:  left translation with $g$ in $M_r(D)$ is an $F$ linear operator $L(g)$ and it has a characteristic 
polynomial $P_{L(g)}$. It 
is a monic polynomial of degree $n^2$. One can  prove that this polynomial is always the power $n$ of a monic polynomial which is, again, $P_g$.

Let $g\in G$, we say $g$  is {\bf elliptic} if $P_g$ is irreducible and has simple roots in $\bar{F}$. We say $g$ is {\bf regular semisimple} if $P_g$ has 
simple roots in 
$\bar{F}$. Let $\tg$ be the set of regular semisimple elements of $G$, which we familiarly call the {\bf regular set}. If $g\in\tg$, then $P_g$ is also the
 minimal polynomial of 
$g$ over $F$. If $g,h\in\tg$, then $h$ is conjugated to $g$ if and only if $P_g=P_h$, as showed in the following lemma. Let $\O_G$ be the set of conjugacy 
classes in $G,$ 
$\tilde{\O}_G$ the set of conjugacy classes of regular semisimple elements and $\tilde{\O}^{ell}_G$ the set of conjugacy classes of elliptic elements.

For $k|n$, let ${\mathcal X}_k\subset F[X]$ be the set of monic polynomials $P$ of degree $n$ with distinct non zero roots in $\bar{F}$ and such that, if $P=\prod_iP_i$ is the 
decomposition of  
$P$ in irreducible factors $k$ divides the degree of each $P_i$.

\begin{lemme}
The map $g\mapsto P_g$ is a bijection from $\tilde{\O}_G$ to ${\mathcal X}_d$ and from $\tilde{\O}^{ell}_G$ to ${\mathcal X}_n.$
\end{lemme}

{\bf Proof.} We prove only the bijection between  $\tilde{\O}_G$ and  ${\mathcal X}_d$; the bijection between  $\tilde{\O}^{ell}_G$ and ${\mathcal X}_n$ being obvious.

 First we show that $g\in\tg$ implies $P_g\in {\mathcal X}_d$. We do it by induction on $r$. Let $r=1$. Then, if $g$ is regular semisimple, 
then $P_g$ is irreducible. Indeed,
 we know that $P_g$ is also the minimal polynomial of $g$, and as $D$ is an integral domain, it has to be irreducible. Now assume $r>1$. If
  $P_g$ is irreducible, the result is 
 clear. Assume $P_g=P_1P_2$ with $P_1$ and $P_2$ non constant.
As it is pointed out in \cite{Lang}, XVII sect. 1, if $D'$ is the opposite algebra to $D$ and we consider the left-$D'$-vector space 
$V:=D'^r$ endowed with the canonical basis, 
then the usual way of associating a matrix to a linear map {\it in the commutative case} yields here a left-$D$-linear isomorphism 
$M_r(D)\simeq End_{D'}V$. If $g\in M_r(D)$ we 
denote $f_g$ the associated $D'$-endomorphism. As $P_g(g)=0$, one has $P_g(f_g)=0$. Now $P_1$ and $P_2$ are mutually prime because 
$P_g$ has simple roots in $\bar{F}$. 
Write $UP_1+VP_2=1$ with $U,V\in F[X]$. It is easy to see that, as in the commutative case, $U(f_g)P_1(f_g)$ and $V(f_g)P_2(f_g)$ 
are associated non zero projectors which 
both commute to $f_g$ (because all the coefficients of the polynomials involved are in $F$), and yield a non trivial decomposition 
of $V=V_1\oplus V_2$ of $V$ into a direct sum of spaces 
stable by $f_g$. Base change implies then that $g$ is conjugated with an element of $M_{r_1}(D)\t M_{r_2}(D)\subset M_r(D)$, $r_1+r_2=r$, $r_1r_2\neq 0$. We then 
apply the 
induction assumption. This proves $P_g\in {\mathcal X}_d$.

We show now that the map $g\mapsto P_g$ is injective. As $g\in \tilde{G}$, the subalgebra $F[g]$ of $A$ generated by 
$g$ is isomorphic to $F[X]/(P_g)$ by sending $g$ to the class of ${X}$. So, if $g$ and $g'$ are such that $P_g=P_g'$, 
there is an isomorphism $i: F[g]\to F[g']$ sending $g$ to $g'$.
Assume first that $P_g$ is irreducible. Then $F[g]$ is a field and as $A$ is a simple algebra, the result follows by 
Skolem-Noether theorem which asserts that $i$ is conjugation with an element of $A$. The general case follows then by induction, as before.

We show the surjectivity. Let first $P$ be an {\it irreducible} monic non constant polynomial over $F$ of degree 
divisible by $d$. Assume $P$ has simple roots in $\bar{F}$. Consider the extension $E:=F[X]/(P)$  of $F$ of degree 
equal to $\deg P$. According to \cite{Pie}, Corollary 13.3, there exists a subfield of $M_{\frac{\deg P}{d}} (D)$ 
isomorphic to $E$. So $M_{\frac{\deg P}{d}} (D)$ contains an element $g$, such that $P_g=P$. Moreover $g$ is 
(invertible and) regular semisimple by definition. Now pick up any element $P$ of ${\mathcal X}_d$ and decompose it 
$P=\prod_iP_i$ in 
irreducible factors. By definition, the degree of each $P_i$ is divisible by $d$. For each $i$, let 
$g_i\in M_{\frac{\deg P_i}{d}} (D)$ such that $P_{g_i}=P_i$. Then let $g\in M_r(D)$ be the element in 
the Levi subgroup $\prod_i GL_{\frac{\deg P_i}{d}}(D)$ whose  blocks are the $g_i$. Then $P_g=P$. This proves 
the surjectivity.\qed
\ \\

If $g\in A$ (resp. $g\in G$), $A_g$ (resp. $G_g$) will be the centralizer of $g$ in $A$ (resp. in $G$). If $g\in\tg$, then $\bar{X}\mapsto g$ is an 
embedding of $F[X]/(P_g)$ in 
$A$ with image $A_g$. 
$G_g$ is a maximal torus of $G$, isomorphic to the group $A_g^\t$ of invertible elements of $A_g$. The set $\tg_g$ of regular semisimple elements of 
$G_g$ is a dense subset of 
$G_g$. In the following we will use the lemma:

\begin{lemme}\label{measpreserving}
Let $g\in\tg$ and fix a Haar measure on $G_g$. Let $i$ be a continuous automorphism of $G_g$ such that, for all $h\in\tg_g$, $i(h)$ has the same 
characteristic polynomial as $h$. Then $i$ is measure preserving. 
\end{lemme}

{\bf Proof.} If $i$ is conjugation by an element of $G$ this comes from the fact that the Weyl group is finite, the Haar measure is unique up to a
 scalar in $\r_+^\t$ and a 
finite subgroup of $\r_+^\t$ is trivial. 

In the general case, as $g$ is regular semisimple, $i(g)$ is conjugated to $g$. So composing $i$ with the 
appropriate conjugation, which is measure preserving,  one may then assume that $i(g)=g$. Now there 
is an open neighborhood $V$ of $g$ in $G_g$ such that all the elements of $V$ are regular semisimple, and not conjugated to each 
other (\cite{HC} for example). The map $g\mapsto P_g$ is so injective on $V$. Then $i^{-1}(V)$ has the same property. If $W:=V\cap i^{-1}(V)$ then
 $W$ is an open neighborhood of $g$, and, as $i$ preserves the characteristic polynomial we have to have $i(h)=h$ for all $h\in W$. So the restriction 
 of $i$ to an open set is identity and $i$ is measure preserving.\qed
   
\ \\

\subsection{Transfer of orbits}

We now  change notation in order to fit to the standard literature in this field:  we set $A':=M_r(D)$, $G':=GL_r(D)$, like before, and $A:=M_n(F)$, $G:=GL_n(F)$. We identify the centers of $G$ and $G'$ by the canonical isomorphism and we call it $Z$. 
If $d$ is a positive  integer dividing $n$, we let $\tg^d$ be the set of elements
 $g\in \tg$ such that  $P_{g}\in {\mathcal X}_d$.

We write $g\lra g'$ and we say that $g$ {\bf corresponds} to $g'$ if $g\in \tg^d$, $g'\in\tg'$ and $P_g=P_{g'}$. 

Because  $\tilde{\O}_{G'}$ is in bijection with ${\mathcal X}_d$ and  $\tilde{\O}_{G}$ is in bijection with 
${\mathcal X}_1$, the inclusion ${\mathcal X}_d\subset {\mathcal X}_1$  induces an injective map from $\tilde{\O}_{G'}$ to $\tilde{\O}_G$ associated to  the 
previous correspondence. 

\subsection{Transfer of centralizers}

On tori of type $G_g$, $g\in\tg$, of $G$ we fix Haar measures such that if two such tori are conjugated then the measures are compatible with the conjugation.
 Moreover, if $G_g/Z$ is compact (i.e. $g$ is elliptic), we assume the measure gives volume one to $G_g/Z$. This is well defined thanks to the lemma 
 \ref{measpreserving}.

We are going to fix Haar measures on tori $G'_{g'}$, $g'\in\tg'$, of $G'$. If $g'\in \tg'$, let $g\in\tg$ such that $g\lra g'$. Then $P_g=P_{g'}$ and 
we get  canonical isomorphisms $A_g\simeq F[X]/(P_g)\simeq A'_{g'}$ which preserve the characteristic polynomial. Then we get an 
isomorphism $G_g \simeq G'_{g'}$ (both are isomorphic to $(F[X]/(P_g))^\t$) and we use this isomorphism to define a Haar measure on $G'_{g'}$ through 
transfer from $G_g$. This is well defined (does not depend of choices) thanks to the lemma \ref{measpreserving}. Moreover, if $G'_{g'}$ and $G'_{h'}$ are
 conjugated then the measures are compatible with the conjugation and if $G'_{g'}/Z$ is compact (i.e. $g'$ is elliptic), the measure gives volume one 
 to $G'_{g'}/Z$.

\subsection{Transfer of functions}

If $C$ is a non empty subset of $G$, we denote

- $1_C$  the characteristic function of $C$, 

- $Ad(G)C$  the set of all elements of $G$ which are conjugated to an element of $C$, 

- $H(C)$ the set of complex functions on $G$ which are locally constant and has compact support included in $C$.\\
We denote $Supp(f)$ the support of a function $f$. 

If $f\in H(G)$, then we define
 the orbital integral of $f$ in a 
point $g\in \tg$ by
$$\Phi(f,g):=\int_{G_{g}\bc G}f(x^{-1}g x)dx$$
where $dx$ is the quotient measure. The integral is convergent (\cite{Lau} proposition (4.8.9)). 
$\Phi(f,\cdot )$ is locally constant on 
$\tg$ and stable by conjugation under $G$. If $f\in H(\tg)$, then we have 
$Supp(\Phi(f,\cdot))\subset Ad(G)Supp(f)\subset \tilde{G}$.

According to the Harish-Chandra submersion theorem \cite{HC}, every $g\in \tg$ has a neighborhood $V$ in $\tg$ such that there is
 an open compact subgroup $K_g$ of $G$ and a neighborhood $V_g$ of $g$ in $G_g\cap \tg$ such that the map $K_g\t V_g\to V$ defined by $(k,x)\mapsto k^{-1}xk$ is an 
 isomorphism. We will call such a neighborhood a HC-neighborhood. Notice that the orbital integral $\Phi(1_V,\cdot)$ of the characteristic function of $V$ is a scalar multiple of $1_{Ad(G)V}$.
A classical application of Harish-Chandra submersion theorem is then the following lemma:

\begin{lemme}\label{HC}
Let $C$ be an open compact subset of $\tilde{G}$. Let $\Phi:G\rightarrow \mathbb{C}$ be a locally constant function stable by conjugation, such 
that $Supp(\Phi)\subset Ad(G)C$. Then $\Phi$ is the orbital integral of a function $f\in H(C)$. 
\end{lemme}  
\ \\
{\bf Proof.} Let $C=\cup_{j\in J}V_j$ be a covering of $C$ with open sets $V_j$ such that every set $V_j$ is included in 
a HC-neighborhood. One may write $C=\coprod_{i=1}^k U_i$ where $U_i$ is open compact, $\Phi$ is constant on $U_i$ and for every 
$i$ there exists $j$ such that $U_i\subset V_j$ (\cite{Re}, Lemma II.1.1.ii). Then the orbital integral 
of $1_{U_i}$ is constant and non zero on $Ad(G)U_i$ and so there is a scalar $\lambda_i$ such that 
the orbital integral $\Phi(\lambda_i 1_{U_i},\cdot )$ is equal to $\Phi$ on $Ad(G)U_i$. The function
 $f:=\sum_{i=1}^k \lambda_i 1_{U_i}$ has the required property.\qed

We adopt the same notation with $G'$ instead of $G$. The same results are true for $G'$. 
We write $f\lra f'$ and say that $f$ {\bf corresponds} to $f'$ if $f\in H(\tg^d)$, $f'\in H(\tg')$, and we have

- $\Phi(f,g)=\Phi(f',g'),  \forall g\in \tg,\forall g'\in \tg, g\lra g'$,

- $\Phi(f,g)=0$ if $g\in \tg\bc\tg^d$.\\
A consequence of the lemma \ref{HC} is the following:

\begin{prop}
(a) If $f\in H(\tg^d)$, then there exists $f'\in H(\tg')$ such that $f\lra f'$.

(b) If $f'\in H(\tg')$ then there exists $f\in H(\tg^d)$ such that $f\lra f'$.
\end{prop}

\subsection{Transfer of unitary representations} If $\pi$ is a smooth irreducible representation and $f\in H(G)$, one defines the finite rank 
 operator $\pi(f)$ by the usual formula $\pi(f):=\int_G f(g)\pi(g)dg$. If $\pi$ and $\pi'$ are isomorphic, then $\tr\pi(f)=\tr\pi'(f)$. 
Let $Irr(G)$ be the set of isomorphy classes of smooth irreducible representations of $G$ and $Irr_u(G)$ the 
subset of unitarizable (classes of) representations. 
Harish-Chandra (\cite{HC}) attached to the smooth irreducible representation $\pi$ its character $\chi_\pi$, defined in {\it any} characteristic, which verifies:

- $\chi_\pi$ is a locally constant  function from $\tg$ to $\cc$, which is stable by conjugation 

- if $f\in H(\tg)$, then for every representation $\sigma$ in the isomorphy class of $\pi$ one has
 $\tr\sigma(f)=\int_{\tg}f(g)\chi_\pi(g)dg$\\
(original result by Harish-Chandra \cite{HC2}; see also \cite{DS} and \cite{MS}).

This holds also for $G'$, and we define $Irr(G')$, $Irr_u(G')$ and $\chi_\pi$ for $\pi\in Irr(G')$ in the same way.

Harish-Chandra (\cite{HC2}) proved, {\it when the characteristic of $F$ is zero}, that $\tr\sigma(f)=\int_{\tg}f(g)\chi_\pi(g)dg$ for any $f\in H(G)$ (see also \cite{DS}) (resp. $f\in H(G')$).
This was also proved to hold for $G$ (\cite{Le1}) and $G'$ (\cite{Bacrelle}, \cite{Le2}) {\it when the characteristic of $F$ is positive}.

We will frequently identify irreducible representations with their class in $Irr(G)$ when using notions which  are invariant
 under isomorphism. Let 
$Irr_u^d(G)$ be the set of representations of $\pi\in Irr_u(G)$ such that the restriction of $\chi_\pi$ to 
$\tg^d$ is not null. We have the following theorem, proved in \cite{BHLS}, which is  a local Jacquet-Langlands 
transfer in 
positive characteristic  for all irreducible unitary representations generalizing \cite{DKV}:

\begin{theo}\label{bhls}
There is a unique map $\lj:Irr_u^d(G)\to Irr_u(G')$ such that, for every $\pi\in Irr_u^d(G)$ there exist $\varepsilon(\pi)\in\{-1,1\}$ such that  
$$\chi_\pi(g)=\varepsilon(\pi) \chi_{\lj(\pi)}(g')$$
for all $g\lra g'$.
\end{theo}
\ \\

In general, the map $\lj$ is neither  injective nor surjective.

\def\gf{G(F)}
\def\ga{G(\aa)}

\section{Main result}

\subsection{Basic facts}

Let $F$ be a global field of characteristic $p$ i.e.   a finite extension of the field
of fractions ${\mathbb F}_p(X)$. Fix an algebraic closure $\bar{F}$ of $F$. 
For each place $v$ of $F$, let $F_v$ be the completion of $F$ at $v$,  $O_v$ be the ring of integers of $F_v$, and fix once for all
 an algebraic closure $\bar{F}_v$ of $F_v$.
 
 Let $D$ a
central division algebra over $F$ of dimension $n^2$. We set $A=D,$  and for each place $v$ of $F$ let $A_v:=D\otimes_{F} F_v$. 
$A_v$ is a
central simple algebra over $F_v$ and by Wedderburn theorem  $A_v\simeq M_{n_v}(D_v)$ for some positive integer
$n_v$ and some central division algebra $D_v$ of dimension $d_v^2$
over $F_v$ such that $n_v d_v=n$. We will fix once and for all an
isomorphism and identify these two algebras.
 We will denote $O'_v$ the ring of integers of $D_v$.

 We say that $D$
{\bf splits} at a place $v$ if $d_v=1$. The set $V$ of places
where $D$ does not split is finite and  it is known by the class field theory that  $n$ is the least common multiple of the
$d_v$ over all the places $v\in V$.

Let $G$ be the group $GL_n(F)$, and for
  each place $v$ of $F,$ let $G_v$ be the group $GL_n(F_v)$. Let then $K_v$ be the maximal compact subgroup $GL_n(O_v)$ of $G_v$. 
  
Let $G'$ be the group $D^{\times}$;  for every place $v\in
V$, set $G'_v= A_v^\times= GL_{r_v}(D_v)$. Set then $K'_v:=GL_{r_v}(O'_v)$ a maximal compact 
subgroup of $G'_v$. For $v\notin V$, we fix once for all an isomorphism $A_v\simeq M_n(F_v)$ and we identify these algebras. Notice that such an isomorphism is, by Skolem-Noether theorem, unique up to a conjugation by an invertible element of the algebra. Identify consequently  $G'_v$ and $G_v$ and set $K'_v:=K_v$.

Let $\aa$ be the ring of adeles of $F$ and  denote  $G(\aa)$  the adelic  group of $G$ with respect to the $K_v$. We consider $G$ 
  as a subgroup of $G(\aa)$ by the diagonal embedding. Let $Z$ be the center of $G$; it is identified with $F^{\times}$, and for each place 
  $v$, let $Z_v$ be the center of $G_v$ also 
  identified with $F_v^{\times}$. Let $Z(\aa)$ be the center of $G(\aa)$, also identified with the adelic  group of 
  $Z$ with respect to open compact subgroups $K_v\cap Z_v$. $Z(\aa)$ identifies with the group of ideles $\aa^{\times}$ of $F.$ 
   For every place $v$ of $F$, fix the Haar measure $dg_v$  on $G_v$ such that the volume of $K_v$ is one, and 
   $dz_v$ on $Z_v$ such that the volume of $Z_v\cap K_v$ is one. On $G(\aa)$ (resp. $Z(\aa)$) consider then the 
   unique product Haar measure $dg$ (resp. $dz$). 

One defines a group morphism $\text{deg}:\aa^{\times}\rightarrow \z$, as defined page 15 of \cite{Laf} or Part II of \cite{Lau} page 3, by
\begin{equation*} deg(a)=\sum_{v}deg(v) v(a_v)
\end{equation*}
where $a=(a_v)_v$, $\kappa_v$ is the residual field of $F_v$, $deg(v)$ denotes 
the dimension of $\kappa_v$ over $\mathbb{F}_p$ and the sum is taken over all places $v$ of $F.$

This morphism is surjective
(Lemma 9.1.4 page 3 of \cite{Lau}). 
We let $a=(a_v)_v\in \aa^{\times}$ an idele of degree $1$. According to the Lemme 1 page 48 of \cite{Laf}, {\it we may assume that $a_v=1$ outside a finite set of places  $T_a$ such that $T_a\cap V=\emptyset$}. This is not essential for the proof, but it highly simplifies computations.
 Let $J:=a^{\z}$ the subgroup of $Z(\aa)\simeq \aa^\t$ generated by $a$. It is not a restricted product over all the places, but may be written as a product $J_{T_a}\times \{1\}$, where $J_{T_a}$ is a subgroup of $\times_{v\in T_a} Z_v$ and $\{1\}$ is to be understood as the trivial subgroup of the restricted product $\t '_{v\notin T_a} Z_v$.

We denote $G'(\aa)$ the adelic group of $G'$ with respect to the $K'_v.$
We consider $G'(F)$ as a subgroup of $G'(\aa)$ by the diagonal embedding.

There are canonical isomorphisms between the center
 of $G$ and the center of $G'$, and, for all place $v$, between the center of $G_v$ and the center of $G'_v,$ so we will identify them. 
  The same is true for the center of $G(\aa)$ and the center of $G'(\aa)$ which will be identified.
For every place $v$ of $F$, fix the Haar measure $dg'_v$  on $G'_v$ such that the volume of $K'_v$ is one. 
On $G'(\aa)$ consider then the product Haar measure $dg'$. For the places $v\notin V$, the identification between $G_v$ and $G'_v$
 is compatible with these choices.\\

For the theory of parabolic subgroups of $G$ we adopt the same conventions and notation as in the local case, which are the conventions of \cite{Lau}
 for example
$\Delta=\{1,\cdots,r-1\}$, and to any subset $I\subset \Delta$ we associate a standard parabolic subgroup $P_I(\aa)$, with Levi
 decomposition $P_I(\aa)=M_I(\aa)N_I(\aa)$ etc..  
If $P=P_I$ is a standard parabolic subgroup of $GL_n(\aa)$, we will sometimes write $M_P:=M_I$ for the Levi component of $P$ and $N_P:=N_I$ for its
 unipotent radical, $P=M_PN_P$. Same notation over $F$: $P_I(F)$ etc.. $M_0:=M_\emptyset$ is the minimal standard Levi subgroup made of 
 diagonal matrices of $GL_n$ and ${\mathcal P}_0$ will denote the finite set of all parabolic subgroups, standard or not, containing $M_0$. Let ${\mathcal P}_0^s$ be the subset of ${\mathcal P}_0$ made of standard parabolic subgroups. Every $P\in {\mathcal P}_0$ has a Levi component which is a standard Levi subgroup, denoted $M_P$. Then, if $r_I=(r_1,...,r_k)$ is the partition associated to $M_P$, we define a homomorphism 
 $deg_{M_P}: M_P(\aa)\rightarrow {\mathbb Z}^k$ by:
 \begin{equation*}
 deg_{M_P}(g)=(deg(det(g_1)),...,deg(det(g_k)))
 \end{equation*}
 where $g=diag(g_1,...,g_k)$ is its  block decomposition, (note:we use here the normalization of L.Lafforgue \cite{Laf} p280).

\subsection{Automorphic representations}

In this subsection we follow \cite{Lau} and \cite{Laf}.
We will be concerned with the representation  of $G(\aa)$  acting on  the space of  functions on 
$G(F)\bc G(\aa)/J$ by right translation. We endow 
$G(F)\bc G(\aa)/J$ and $G'(F)\bc G'(\aa)/J$ with the quotient measures. According to \cite{Laf}, III.6.Lemme 5, $G'(F)\bc G'(\aa)/J$ is compact,
 $G(F)\bc G(\aa)/J$ has finite measure, and they both have the same measure.
We denote $R_{G}$ the representation of $G(\aa)$ acting on the space $L^2(G(F)\bc G(\aa)/J)$ by right translations.

We have a variant,  for any  parabolic subgroup $P\in {\mathcal P}_0$ of $G,$ of this representation  which is the representation of  $G(\aa)$ acting on the space
of $L^2(M_P(F)N_P(\aa)\bc G(\aa)/J)$ by right translation. We denote $R_{G}^{P}$ this representation. Note
that $R_G^G=R_G.$

We   also need the representations of $M_P(\aa )$ on the space of $L^2(M_P(F)\bc M_P(\aa)/J)$ which is important for
defining the  notion of discrete pair.

Let $P=MN$ be a  parabolic subgroup of $G$, $Z_M$ the center of $M$ and  $\chi:Z_M(F)\bc Z_M(\aa)/J\rightarrow \cc$ a 
 smooth character, let $K'\subset K$ an open subgroup of $K$, we denote $L_{K'}^2(M(F)\bc M(\aa)/J,\chi)$ the space of (necessarily locally constant) 
 functions $f$ on $M(\aa)/J$ with values in $\cc$ such
 that:
\begin{itemize}
\item $f(zm)=\chi(z)f(m), \forall z\in Z_M(\aa)/J, \forall m\in M(\aa)/J$
\item $f$ is  invariant on the left by $M(F)$
\item $f$ is invariant  on the right by $K'\cap M_P(\aa)$ 
\item $f$ is of finite norm in the sense of Lafforgue \cite{Laf} page 282.
\end{itemize}

We denote $L^2_{\infty}(M(F)\bc M(\aa)/J,\chi)$ the inductive limit $
\lim\limits_{\rightarrow}L^2_{K'}(M(F)\bc M(\aa)/J,\chi).$ 

One therefore obtains a representation $R_{M_P, \chi}$ of $M_P(\aa )$  
acting on $L^2_{\infty}(M(F)\bc M(\aa)/J,\chi)$ by right translation.

An irreducible  subrepresentation  of  $R_{M_P, \chi}$  is called a {\bf discrete series} of $M(\aa)$. The subspace of $R_{M_P, \chi}$ generated by
 all irreducible subrepresentations is denoted $L^2_{\infty}(M(F)\bc M(\aa)/J,\chi)_{disc}.$ 
 The  isotypical components  of $L^2_{\infty}(M(F)\bc M(\aa)/J,\chi)_{disc}$ are called the {\bf discrete components}.

A {\bf discrete pair} is a couple $(P,\pi)$ where $P\in {\mathcal P}_0$ and $\pi$ a discrete component of  $R_{M_P, \chi_\pi}$ for
some 
central character $\chi_{\pi}$ (which is necessarily the central character of $\pi$).

Every discrete series 
$\pi$ of $G(\aa)$ is isomorphic with a restricted Hilbertian tensor
product of (smooth) irreducible unitary representations $\pi_v$ of
the groups $G_v$ as explained   in \cite{Fl1}. Each representation $\pi_v$ is
determined by $\pi$ up to isomorphism and is called the {\bf local
component of $\pi$ at the place $v$}. For almost all places
$v$, $\pi_v$ has a non zero fixed vector under $K_v$. We say then
$\pi_v$ is {\bf spherical}. 

The same definitions and properties hold for $M_P(\aa)$ and for $G'(\aa)$.\\

\subsection{Relation with the classical setting}\label{relations}

This setting is slightly different from the classical one (\cite{MW} or \cite{Ar}) and references therein. This is because the quotient with
 this subgroup $J$ is very convenient in non zero characteristic. As the cornerstone of our proof is the Theorem 12, VI.2 from \cite{Laf}, we need this definition. Let us
  explain quickly the link with the classical setting: let us say that a discrete series in the sense of \cite{MW} is {\it cctJ} if it has 
  {\it central character trivial on $J$}. Then the discrete series of $G(\aa)$ as defined here correspond exactly to the cctJ discrete series in the 
  classical setting. In particular, the multiplicity one theorem holds for $G(\aa)$ and our discrete series (not {\it a priori} for $G'(\aa)$ but we 
  will prove it here). The other way round, a discrete series in the classical setting is always a twist with a character of a cctJ (following lemma) so 
  proving the Jacquet-Langlands correspondence in Lafforgue's setting leads also to the desired result in the classical setting.

\begin{lemme}\label{change}
{\rm (a)} Let $\chi$ be a character of $F^\t\bc Z(\aa)$. Then there exists a character $\chi'$ of $G(F)\bc G(\aa)$ (or $G(F)\bc G'(\aa)$), 
$\chi'=p^{s\deg\circ\det}$ where $s$ is a complex number, such that $\chi\chi'$ is trivial on $J$.

{\rm (b)} If $\pi$ is a discrete series in the sense of \cite{MW} with central character $\chi$, if $\chi'$ is like in (a), then $\chi'^{-1}\otimes \pi$ is cctJ.
\end{lemme}

{\bf Proof.} (a) We search for $s$ such that $\chi'(a)\chi(a)=1$. Recall $\deg\det(a)=n$. Let $z$ be a $n$-th root of $\chi(a)$. It is enough to chose $s\in \cc$ such that $p^{s}=z^{-1}$. (b) is obvious.\qed
\ \\

\subsection{Claim of the correspondence}

If $\pi$ is a discrete series of $G(\aa)$, we say $\pi$ is $D$-{\bf compatible} if the local components of $\pi$ verify : for all $v\in V$, $\pi_v\in Irr_u^{d_v}(G_v)$.

We will prove the following theorems:

\def\f{{\bf G}}
\begin{theo}\label{correspondence} {\bf Global Jacquet-Langlands correspondence.}

There exists a unique map $\f$ from the set of $D$-compatible discrete series 
of $G(\aa)$ to the set of discrete series of $G'(\aa)$ such that for all discrete series $\pi$ of 
$G(\aa)$ if $\pi'=\f(\pi)$ then 

- $\lj_v(\pi_v)=\pi'_v$ for all places $v\in V$, and 

- $\pi_v=\pi'_v$ for all places $v\notin V$\\ 
where  $\lj_v$ denote the local Langlands-Jacquet correspondence at place $v$ of theorem \ref{bhls}.\\ 
The map $\f$ is bijective. 
\end{theo}

\begin{theo}\label{multiplicity one } {\bf Multiplicity one Theorems for $G'(\aa)$.}

{\rm (a)} If $\pi'$ is a discrete series of $G'(\aa)$, then $\pi'$ appears with multiplicity 
one in the discrete spectrum (multiplicity one theorem).

{\rm (b)} If $\pi'$, $\pi''$ are discrete series of $G'(\aa)$ such that $\pi'_v\simeq \pi''_v$ for 
almost all place $v$, then $\pi'=\pi''$ as subrepresentations of $L^2(G'(F)\bc G'(\aa)/J)$ (strong multiplicity one theorem).
\end{theo}

The rest of the paper is devoted to the proof of these theorems. We will work with the Laumon-Lafforgue trace formula. Then, the lemma 
\ref{change} (b) allows one to transpose the theorem in the classical setting.

\section{The proof}

\subsection{Transfer of elliptic global orbits} 

Characteristic polynomials are defined in the global case like in the local case. \cite{Pie} does not 
treat explicitly the global characteristic $p$ case, but the reader may find it in \cite{We}. If $g\in D$ has characteristic polynomial
 $P_g\in F[X]\subset F_v[X]$, then $P_g$ is the characteristic polynomial of $g$ as an element of  $A_v$ for all $v$ since an embedding $D\hookrightarrow M_k(F)$ uniquely extends to a continuous embedding $A_v\hookrightarrow M_k(F_v)$.

We say an element $g\in G(F)$ (resp. $g\in G'(F)$) is {\bf elliptic} if the characteristic polynomial of $g$ is irreducible over $F$ and has simple
 roots in $\bar{F}$. Let $\tilde{O}_{G(F)}^{ell}$ (resp. $\tilde{O}_{G'(F)}^{ell}$) be the set of elliptic conjugacy classes in $G(F)$ (resp. $G'(F)$). 

Let ${\mathbb X}$ be the set of monic polynomials $P$ of degree $n$ with coefficients in $F$ such that $P$ is irreducible over $F$ and has 
simple roots in $\bar{F}$. Let ${\mathbb X}_D$ be the subset of polynomials $P\in {\mathbb X}$ such that for all place
 $v\in V$, if $P=\prod P_i$ is the decomposition in irreducible factors of $P$ over $F_v$, then for all $i$ the integer $d_v$ divides $\deg P_i$.

Then we have

\begin{lemme}\label{ellipticglobal}
(a) The map $g\mapsto P_g$ induces a bijection from $\tilde{O}_{G(F)}^{ell}$ to ${\mathbb X}$.

(b) The map $g\mapsto P_g$ induces a bijection from $\tilde{O}_{G'(F)}^{ell}$ to ${\mathbb X}_D$.
\end{lemme} 

{\bf Proof.}
(b) 
The fact that the map $\tilde{O}_{G'(F)}\to {\mathbb X}_D$ is well defined, i.e. takes values in ${\mathbb X}_D$, comes from the local analogous
 result. 

The map is injective (this may be proved as in the local case, using the Skolem-Noether theorem). 

The map is surjective: it is consequence of a result of class field theory: Let $P\in {\mathbb X}_D$ and set $L:= F[X]/(P)$ which we see as
 an extension of $F$. Then $L\otimes F_v$ is a product of fields, isomorphic to $F_v[X]/(P_i)$, where $P_i$ are the prime factors 
 of $P$ over $F_v$. The condition $P\in {\mathbb X}_D$ implies that the extension $L/F$ verifies the condition (ii) of Proposition 5, 
 \cite{We} XIII sect. 3, page 253. The equivalence stated in this proposition between (ii) and (iii) implies that $L$ is isomorphic 
 to a subfield of $D$. The element $\bar{X}\in F[X]/(P)=L$ is then sent to an element $g\in D$ whose characteristic polynomial 
 is $P$, as required. This proves the map is surjective.

(a) is now a particular case of (b). The surjectivity in (a), however, is easier to prove using the companion matrix.\qed
\ \\

Let $\tilde{G}(F)^D$ be the set of $g\in G(F)$ such that $P_g\in {\mathbb X}_D$. Let $\tilde{O}_{G(F)}^D$ be the set of conjugacy classes of 
$\tilde{G}(F)^D$.\\
\ \\

\subsection{Transfer of global functions}

Let
$H(G(\aa))$ (resp. $H(G'(\aa)$) be the set of functions $f:G(\aa)\to\cc$ (resp. $f:G'(\aa)\to\cc$) such that $f$
is a product $f=\otimes_v f_v$  over all places $v$ of $F$, where $f_v\in H(G_v)$ (resp. $f_v\in H(G'_v)$) for all $v$, and,
 for almost all
$v$, $f_v$ is the characteristic function of $K_v$ (resp. of $K'_v$). We write then
$f=(f_v)_v$.

Let $\tilde{G}(\aa)^D$ be the set of elements $g\in G(\aa)$ such that for every place $v\in V$ we have
 $g_v\in \tilde{G}_v^{d_v}$, which is also the set of elements $g\in G(\aa)$ such that, for all place 
 $v\in V$, if $P_{g_v}=\prod P_i$ is the decomposition of the characteristic polynomial of $g_v$ in a product of 
 irreducible polynomials over $F_v$, then $d_v$ divides the degree of every $P_i$.

Let $H(\tg(\aa)^D)$ be the subset of $H(G(\aa))$ made of functions $f=(f_v)_v\in H(G(\aa))$ such that, 
for all $v\in V$, $f_v\in H(\tg_v^{d_v})$. Let $H(\tg'(\aa))$ be the set of functions $f'=(f'_v)_v\in H(G'(\aa))$ 
such that, for all $v\in V$, $f_v\in H(\tg'_v)$. 

If $f=(f_v)_v\in H(G(\aa))$ and $f'\in H(G'(\aa))$, we write 
$f \stackrel{\aa}{\lra} f'$  if $f\in H(\tg(\aa)^D)$, $f'\in H(\tg'(\aa))$ and, for all $v\in V$, $f_v\lra f'_v$, 
and for all $v\notin V$, $f_v=f'_v$. If $f\in H(\tg(\aa)^D)$, then there exists $f'\in H(\tg'(\aa))$ such that $f \stackrel{\aa}{\lra} f'$ and if 
$f'\in H(\tg'(\aa))$, then there exists $f\in H(\tg(\aa)^D)$ such that $f \stackrel{\aa}{\lra} f'$. This is a direct consequence of the local transfer of functions.

\begin{prop}\label{poly} If $f \in H(\tg(\aa)^D)$, we have

(a) If $g\in G(\aa)\ \bc\ \tilde{G}(\aa)^D$, then $f(g)=0$.

(b) If $g\in G(\aa)$ is conjugated to an element of a standard proper parabolic subgroup $P(\aa)$ of $G(\aa)$, then $f(g)=0$.

(c) $G(F)\cap \tilde{G}(\aa)^D=\tilde{G}(F)^D$.
\end{prop}

{\bf Proof.} (a) and (c) are obvious.

(b) Assume $P(\aa)$ be the standard proper parabolic subgroup of $G(\aa)$ associated to the partition $(r_1,r_2,...,r_k)$.
If $g\in G(\aa)$ is conjugated to an element of $P(\aa)$, then, for every place $v$, the characteristic polynomial 
of $g$ breaks into a product of polynomials of degrees $r_1,r_2,...,r_k$. But there exists a place $v_0\in V$ and
 there exists $i$ such that  $d_{v_0}$ does not divide  $r_i$ (by class field theory, the least common multiple 
 of all $d_v$ is $n$). So $f(g)=0$ by (a).\qed

For every $o\in \tilde{O}_{G(F)}^D$ fix once and for all an element $\gamma_o\in o$.
Let $G_{\gamma_o}$ denote the centralizer of $\gamma_o$ in $G$. 
The centralizer $G_{\gamma_o}(\aa)$ of $\gamma_o$ in $G(\aa)$ is the restricted product of the local centralizers $G_{v,\gamma_o}$. 
These local tori are endowed with measures like in the previous section, and $G_{\gamma_o}(\aa)$ is endowed with
the product measure. 
For $f\in H(\tilde{G}(\aa))$ and $o\in O_{G(F)}^D$ we consider  the orbital integral
$$\Phi(f;\gamma_o)=\int_{G_{\gamma_o}(\aa)\bc\ga} f(x^{-1}\gamma_ox)\, dx$$
where $dx$ is the quotient measure. Then $\Phi(f;\gamma_o)$ is the product of local orbital integral $\Phi(f_v;\gamma_o)$.
 We will also have to use orbital integrals $\Phi(f;z\gamma_o)$, where $z\in J$. As $J\subset Z(\aa)$, we 
 have $G_{\gamma_o}(\aa)=G_{z\gamma_o}(\aa)$.

For every $o'\in \tilde{O}_{G'(F)}$ fix once and for all an element $\gamma_{o'}\in o'$. For $f'\in H(\tilde{G'}(\aa))$ and 
$\gamma_{o'}$ we define the orbital integral $\Phi(f';\gamma_{o'})$ in the same way.

\subsection{Trace formula in characteristic $p$}

Laumon and Lafforgue developed, following ideas of Arthur, a trace formula in non zero characteristic.
In this section we review the trace formula for $G(\aa)$ in characteristic $p$ from \cite{Laf} (our Theorem \ref{lafforgue}). This section is devoted to the definition of the ingredients of the formula (we show in the next subsection that  most of them are null for suitable functions).

Let $h:G(\aa)/J\to \cc$ be  a  locally constant function with compact support and $P$ a  parabolic subgroup of $G$.  The convolution
operator $\varphi\mapsto \varphi\ast h$ in the space of square  integrable complex functions on
$M_P(F)N_P(\aa)\bc G(\aa)/J $  is the operator $R_G^P(\check{h})$ where $\check{h}(g)=h(g^{-1})$. It is 
an integral  operator with kernel given by:
\begin{equation*}
K_{h,P}(x,y)=\int_{N_P(\aa)}\sum_{\gamma\in M_P(F)} h(y^{-1}\gamma n_P x)dn_P.
\end{equation*}
We set $K_h:=K_{h,G}.$

Because the function $x\mapsto K_h(x,x)$ is not integrable in general, Arthur defined a notion of 
truncated trace as follows.

We define ${\mathfrak a}_{\emptyset}={\mathbb Q}^r$ and  for $i=1,...,r-1$ let $\alpha_i$ be the linear form on
 ${\mathfrak a}_{\emptyset}$ defined by $\alpha_i(h)=h_i-h_{i+1}.$

If $I\subset \Delta$ we denote 
${\mathfrak a}_{I}=\{x\in {\mathfrak a}_{\emptyset}, \alpha_i(x)=0, \forall i\in I\}$, so
${\mathfrak a}_{\Delta}={\mathbb Q}(1,...,1).$
Let $r_I$ be the partition associated to $I$. The   projection
 ${\mathfrak a}_{I}\rightarrow {\mathbb Q}^k, h\mapsto (h_{r_1}, h_{r_1+r_2},...,h_{r_1+r_2+...+r_k})$ is an isomorphism.
We denote $\lambda_i$  for $i=1,\cdots,r-1,$the fundamental weights, linear forms on ${\mathfrak a}_{\emptyset}$, vanishing on 
${\mathfrak a}_{\Delta}$ and defined by 
$\lambda_i(h)=h_1+...+h_i-\frac{i}{r}(h_1+...+h_r).$
 
For each $J\subset I\subset \Delta$ we have ${\mathfrak a}_J=
{\mathfrak a}_I\oplus {\mathfrak a}_J^I$ where
 ${\mathfrak a}_J^I=\{h\in {\mathfrak a}_J, h_1+...+h_r=0, \lambda_i(h)=0, \forall i\in \Delta\setminus I\}.$

Arthur defines a   function  $\hat{\tau}_J^I$ from ${\mathfrak a}_{\emptyset}$ to $\{0,1\}$ characteristic function of the cone

\begin{equation*}
{\mathfrak a}_I+\{h\in{\mathfrak a}_J^I, \lambda_i(h)>0, \forall i\in I\setminus J\}+{\mathfrak a}_{\emptyset}^J.\end{equation*}

If $g\in G(\aa)$ we can write     $g=n_{\emptyset}m_{\emptyset}k$ with
$ n_{\emptyset}\in N_{\emptyset},  m_{\emptyset}\in M_{\emptyset}$ and $k\in K.$
$ m_{\emptyset}$ is uniquely defined up to multiplication on the right by element of $M_{\emptyset}\cap K.$

Therefore one can define a map $H_{\emptyset}:G(\aa)\rightarrow {\mathfrak a}_{\emptyset},$ with $H_{\emptyset}(g)=
deg_{M_{\emptyset}}(m_{\emptyset}).$

Let $T\in {\mathfrak a}_{\emptyset},$ one defines the Arthur truncated diagonal kernel   
 as being the function on $G(\aa)$ defined by:
\begin{equation*}\label{kernel}
K_h^T(x,x)= \sum_{P\in { \mathcal P}_0^s} (-1)^{\vert P\vert -1}\sum_{\delta\in P(F) \bc G(F)} 
K_{h,P}(\delta x, \delta x) {\bf 1}_{P}^T(\delta x)
\end{equation*}
where $x\in G(\aa)$,  $ {\bf 1}_{P}^T$ are the functions on $G(\aa)$ defined by
 $ {\bf 1}_{P}^T(g)=\hat{\tau}_I^{\Delta}(H_{\emptyset}(g)-T)$ 
with $P=P_I, I \subset \Delta. $
This is well defined because for fixed $x$, the sum over $\delta$ is finite (For characteristic zero this was proved by Arthur \cite{Ar3} Lemma 5.1, in positive characteristic it is the lemma 11.1.1 of \cite{Lau}).

The Arthur truncated diagonal kernel is a compactly  supported function  on 
$G(F) \bc G(\aa)/J$  according to the Proposition 11 page 227 of \cite{Laf}. Therefore 
 one can define the truncated trace of $R_G(\check{h})$ denoted $Tr^{T}(h)$ as being 
\begin{equation*}\label{trace}
Tr^{T}(h)=\int_{G(F) \bc G(\aa)/J}  K_h^T(x,x) dx.
\end{equation*}
This is denoted $Tr^{\leq p}(h)$ in Lafforgue  \cite{Laf} where $p$ is a polygon defined by $T$ on page 221.

We now recall the results on the spectral side for general $h.$  

We need some definitions.

Let $P\in {\mathcal P}_0$, one denotes $\Lambda_P$ the abelian complex Lie group 
(of dimension $\vert P \vert-1$) of complex characters 
$\chi:M_P(\aa)/J\rightarrow \cc^\times$ which factorize through the surjective homomorphism
 $deg_{M_P}:M_P(\aa)/J\rightarrow \z^{\vert P \vert}/(r_1,..., r_{\vert P\vert})\z$.  

As a result each $\lambda_P\in \Lambda_P$ can be written uniquely as
 $\lambda_P(m)=q^{\sum_{j=1}^k \rho_j \frac{deg(m_j)}{r_j}}$ with 
 $(\rho_j)\in \times_{j=1}^k{\mathbb C}/\frac{2i\pi}{r_j log q}{\mathbb Z}$ 
 and $\sum_j \frac{\rho_j}{r_j}\in {\mathbb Z}$.

We have $\Lambda_P= Im\Lambda_P\times Re\Lambda_P$ where $Im\Lambda_P$ (resp . $Re\Lambda_P$) denotes the Lie 
group of unitary characters (resp. of real positive characters).  $Im\Lambda_P$ is a compact group, we denote $d\lambda_P$ its 
normalized 
Haar measure.
If $\lambda_P\in \Lambda_P$ we can canonically extend $\lambda_P$ to a function   on $P(\aa)$ right $N(\aa)$ -invariant and then to a function on $G(\aa)$ right $K-$
invariant using 
the decomposition 
$G(\aa)=M(\aa)N(\aa)K$ (\cite{Laf} p280).

If $(P,\pi)$ is a discrete pair, $K'\subset K$ an open subgroup, we denote

 $L_{K'}^2(M_P(F)N_P(\aa)\bc G(\aa)/J, \pi)$ the space of functions
 $\varphi:M_P(F)N_P(\aa)\bc G(\aa)/J\rightarrow \cc$ right invariant by $K'$ and such that $\forall k\in K$ the function 
 $\varphi_k:M_P(F)\bc M_P(\aa)/J\rightarrow \cc, m\mapsto \rho_{P}(m)^{-1}\varphi(mk)$ 
belongs to the isotypical component $\pi\subset L_{\infty}^2(M_P(F)\bc M_P(\aa)/J, \chi_{\pi}),$
where we have denoted $\rho_{P}$ the square root of the modular character of the group $P(\aa)$.

We let $L_{\infty}^2(M_P(F)N_P(\aa)\bc G(\aa)/J, \pi)$  be the inductive limit 

$\lim\limits_{\rightarrow}L_{K'}^2(M_P(F)N_P(\aa)\bc G(\aa)/J, \pi).$

 One may endow $L_{\infty}^2(M_P(F)N_P(\aa)\bc G(\aa)/J, \pi)$ with a structure of  $G(\aa)$ representation defined by 
 $I_P(\pi)=ind_{P(\aa)}^{G(\aa)}(\pi\otimes 1_{N_P(\aa)} )$.
 
At this point it is convenient to use the notation of Laumon $\pi(\lambda_P)=\pi\otimes \lambda_P,$ when $(P,\pi)$ is a discrete pair and $\lambda_P\in \Lambda_P.$
Because it is sometimes convenient to represent $I_P(\pi(\lambda_P))$ in a vector space independent of $\lambda_P$ one is led to define  the multiplication operator 
 $$[\lambda_P]:
 L_{\infty}^2(M_P(F)N_P(\aa)\bc G(\aa)/J, \pi)\rightarrow L_{\infty}^2(M_P(F)N_P(\aa)\bc G(\aa)/J, \pi(\lambda_P)),$$
  $$f\mapsto f\lambda_P,$$ 
 which is a vector
 space isomorphism (here $\lambda_P$ is the function defined on whole $G(\aa)$ as explained).

 We denote $W$ the Weyl group of $GL_r(F)$, it is isomorphic  to the permutation group $\mathfrak{S}_r,$ and we fix 
 an inclusion $W\subset GL_r(F)$ associating to each permutation the permutation matrix.
If $M$ is a Levi subgroup containing $M_{0}$,  we denote  $W_M=W\cap M.$
If $M, M'$ are two Levi subgroups containing $M_{0}$ we denote $Hom(M,M')$ the set of
 $\sigma\in W_{M'}\bc W/ W_{M}$ such that $\sigma M \sigma^{-1}\subset M'.$ 
If $P,P'$ are two parabolic subgroups element of 
$\mathcal P_{0},$ we denote $Hom(P,P')=Hom(M_P,M_{P'}).$ 
 
Let $(P,\pi)$ a discrete pair and $\sigma:P\rightarrow P'$ an isomorphism,
 each such $\sigma$ is represented by an element $w\in W_{M_P'}\bc W/ W_{M_P}$, and to the representation
 $\pi$ of $M_P(\aa)$ one can associate a representation $\sigma(\pi)$ of 
$M_{P'}(\aa)$ acting on the space $\{\varphi(w^{-1} . w), \varphi \in \pi\}.$
Two discrete pairs $(P,\pi)$ and $(P',\pi')$ are said to be {\bf equivalent} if there exists an isomorphism $\sigma:P\rightarrow P'$ and a character $\lambda_P\in \Lambda_P$ 
such that $\pi'=\sigma(\pi\otimes\lambda_P)$.

Let   $P,P'\in {\mathcal P_0}$ satisfying the condition $M_P=M_{P'}$, and let
 $\varphi\in L_{\infty}^2(M_P(F)N_P(\aa)\bc G(\aa)/J, \pi)$. One defines the function $M_P^{P'}(\varphi,\lambda_P)$ of 
   $g\in  G(\aa) $  a usual by:
\begin{equation*}
M_P^{P'}(\varphi,\lambda_P)(g)=\lambda_{P'}(g)^{-1}
\int_{N_P(\aa)\cap N_{P'}(\aa)\bc N_{P'}(\aa)} \frac{dn_{P'}}{dn_{P,P'}}\varphi(n_{P'}g)\lambda_P(n_{P'}g)
\end{equation*}
where we have denoted $dn_{P,P'}$ the normalized Haar measure on $N_P(\aa)\cap N_{P'}(\aa)$ and  $\frac{dn_{P'}}{dn_{P,P'}}$ the quotient 
measure 
on
 $N_P(\aa)\cap N_{P'}(\aa)\bc N_{P'}(\aa),$ and $\lambda_{P'}$ is the function defined on $G(\aa)$ extending the character  on $M_{P'}(\aa)= M_{P}(\aa)$ defined 
by  $\lambda_{P}$, for the precise definitions see \cite{Laf}.

The integral is convergent under some conditions on $\lambda_P$ recalled in \cite{Laf} page 285 and for 
fixed  $\varphi$, the function  $\lambda_P\mapsto M_P^{P'}(\varphi,\lambda_P)$ admits a meromorphic continuation to 
the whole $\Lambda_P.$

If $\lambda_P$ is such that  $M_P^{P'}(\varphi,\lambda_P)$ is well defined, the function 
$g\mapsto M_P^{P'}(\varphi,\lambda_P)(g)$ belongs to $L^2(M_{P'}(F)N_{P'}(\aa)\bc G(\aa)/J, \pi')$ where $(P',\pi')$ is 
the discrete pair defined by  $\pi'=\sigma(\pi)$  with $\sigma$ associated to an element $w$ of the Weyl group satisfying 
 $M_{P'}=w M_Pw^{-1}=M_P.$ 

The map $[\lambda_{P'}]\circ M_P^{P'}(.,\lambda_P)\circ [\lambda_P]^{-1}:$
\begin{equation*}
 L^2_{\infty}(M_P(F)N_P(\aa)\bc G(\aa)/J, \pi(\lambda_P))\rightarrow 
 L^2_{\infty}(M_{P'}(F)N_{P'}(\aa)\bc G(\aa)/J, \pi'(\lambda_{P'}))
\end{equation*} is an intertwining operator between the representations 
$I_P(\pi\otimes \lambda_{P})$ and $I_{P'}(\pi'\otimes \lambda_{P'}).$

One defines $Fix(P,\pi)$ to be the finite  set of couples $(\tau,\mu_P)$ where $\tau$ is 
an isomorphism $\tau:P\rightarrow P$ and $\mu_P\in \Lambda_P$ such that $\pi$ is isomorphic to 
$\tau(\pi\otimes \mu_P),$
$\mu_P$ is necessarily unitary.
 $Fix(P,\pi)$ can be endowed with a structure of finite group
 (\cite{Laf} page 283) defined as follows:
 $$(\tau ',\mu_P') (\tau,\mu_P) =(\tau '  \tau, \tau^{-1}(\mu_P')\mu_P),
  $$
  and for each $(\tau, \lambda)\in Fix(P,\pi),$ one denotes  $Fix(P,\pi,\tau,\lambda)$ the 
subgroup of elements of  $Fix(P,\pi)$ commuting with $(\tau, \lambda).$
Lafforgue defines a  discrete quadruplet $(P,\pi,\sigma,\lambda_\pi)$ as being a discrete pair
 $(P,\pi)$ and a couple $(\sigma,\lambda_\pi)\in Fix(P,\pi).$
If  $\sigma: P\rightarrow P'$ is an isomorphism, Lafforgue defines  a generalization of the previous intertwining operator \cite{Laf} page 286, 
 $M_{P,\sigma}^{P'}(.,\lambda_P) : 
L^2_{\infty}(M_P(F)N_P(\aa)\bc G(\aa)/J,\pi)\rightarrow L^2_{\infty}(M_{P'}(F)N_{P'}(\aa)\bc G(\aa)/J,\sigma(\pi))$ 
and the operator $[\sigma(\lambda_{P})]\circ M_{P,\sigma}^{\sigma(P)}(.,\lambda_P)\circ  [\lambda_P]^{-1}$ is an intertwining operator between the representation 
$I_P(\pi\otimes \lambda_P)$ and $I_{\sigma(P)}(\sigma(\pi\otimes \lambda_{P})).$

In the following, if $\phi\in L^2_{\infty}(M_P(F)N_P(\aa)\bc G(\aa)/J,\pi)$, 
 $h(\phi,.)$ denotes the analytical function $\lambda_P\mapsto h(\phi,\lambda_P)=((\phi \lambda_P)\star h) \lambda_P^{-1}.$

In \cite{Laf} lemma 6 page 303, Lafforgue  introduces the functions ${\hat{{\bf 1}}}_{P}^T$, which are rational functions on $\Lambda_P$ and satisfy, under
 some condition on $ \lambda_P^0\in Re \Lambda_P$:
\begin{equation*}
(-1)^{\vert P\vert-1}{{\bf 1}}_{P}^T(g)=\int_{Im \Lambda_P} {\hat{{\bf 1}}}_{P}^T(\lambda_P \lambda_P^0)
(\lambda_P \lambda_P^0)(g)d \lambda_P, \forall g\in (M_P(F)N_P(\aa))\bc G(\aa)/J.
\end {equation*}.

He associates (page 299)  to each permutation $\tau$ of $\mathfrak{S}_l$ two surjective maps $\tau^+$ 
(resp.$\tau^-$) from the set $\{1,..,l\}$ to $\{1,...,l^+\}$ (resp. $\{1,...,l^-\}$). 
Lafforgue defines (Lemma 5)  a generalization of  the functions ${{\bf 1}}_{P}^T$, 
denoted ${{\bf 1}}_{P,\tau}^T={{\bf 1}}_{\tau^-(P)}^T (1-{{\bf 1}}_{\tau^+(P)}^T)$ with $\tau\in  \mathfrak{S}_{\vert P\vert}$, 
and their Fourier transform ${\hat{{\bf 1}}}_{P,\tau}^T$ 
which are rational functions on $\Lambda_P$ satisfying the following equality on functions on $M_{\tau(P)}(F)N_{\tau(P)}(\aa)\bc G(\aa)/J$:
\begin{equation*}
\int_{Im \Lambda_P}d\mu_P {\hat{{\bf 1}}}_{P,\tau}^T(\mu_0 \mu_P)
\tau(\mu_0 \mu_P)(.)=(-1)^{\vert\tau^-(P) \vert-1 }
{{\bf 1}}_{P,\tau}^T
\end{equation*}
where $\mu_0\in Re \Lambda_P$ is sufficiently small in the sense of Lafforgue \cite{Laf} page 301.

Finally one obtains the theorem (theorem 12 page 309), where we have used the formula of the Th I.9 contained in \cite{Laf2} which corrects
 two minor misprints (the absence of $\vert \sigma\vert$ 
the incorrect $\tau\sigma(\lambda_\pi^\sigma)$ instead of $\tau\sigma(\lambda_\pi).$).
There is an additional misprint concerning the place of $[\tau\sigma(\lambda_\pi)]$ which should be located on the left.

\begin{theo}\label{lafforgue} (Lafforgue) We have

\begin{equation*}
Tr^{T}(h)=\sum_{(P,\pi,\sigma,\lambda_\pi)}Tr^{T}_{(P,\pi,\sigma,\lambda_{\pi})}(h)
\end{equation*}
where the sum is taken over all good representative of equivalence classes of discrete quadruplet with $\pi$ unitary and 
\begin{eqnarray*}Tr^{T}_{(P,\pi,\sigma,\lambda_\pi)}(h)=
&& \frac{1}{\vert {\text Fix} (P,\pi,\sigma,\lambda_\pi)\vert .\vert \sigma\vert }
\int_{Im_{\Lambda_{P_\sigma}}} d\lambda_\sigma \\
&&\sum_{\lambda_\pi^\sigma} Tr_{L^2(M_P(F)N_P(\aa)\bc G(\aa)/J,\pi)} (M_{(P,\pi,\sigma,\lambda_\pi,\lambda_\pi^\sigma)}^T(.,\lambda_\sigma, h))\label{spectralside}
\end{eqnarray*}
where 
$M_{(P,\pi,\sigma,\lambda_\pi,\lambda_\pi^\sigma)}^T(.,\lambda, h)$ is a finite rank  endomorphism of $L^2(M_P(F)N_P(\aa)\bc G(\aa)/J,\pi)$ defined by 

$$
M_{(P,\pi,\sigma,\lambda_\pi,\lambda_\pi^\sigma)}^T(.,\lambda, h)=
\lim_{\mu_\sigma\rightarrow 1\atop \mu_\sigma\in \Lambda_{P_\sigma}}
\sum_{\tau\in \mathfrak{S}_{\vert P_{\sigma}\vert}}
{\hat{{\bf 1}}}_{P_\sigma,\tau}^T(\mu_\sigma \sigma(\lambda_\pi^\sigma)/\lambda_\pi^\sigma\sigma(\lambda_\pi)) $$
$$([\tau\sigma(\lambda_\pi)]\circ M_{P,\tau\sigma}^{\tau(P)}(.,\lambda\lambda_\pi^\sigma))^{-1}\circ M_{P,\tau}^{\tau(P)}{(.,\lambda_{\pi}^\sigma \lambda/\mu_{\sigma})
\circ h(.,\lambda_\pi^{\sigma}\lambda/\mu_\sigma)} 
.$$
\end{theo}

In this formula we need to explain the notations
 $P_\sigma, \lambda_\pi^\sigma.$

To $(\sigma, \lambda_\pi)\in  Fix(P,\pi)$,
  one associates a parabolic subgroup 
$P_\sigma$ ( \cite{Laf} page 305). $\vert P_\sigma\vert$ is the number of cycles in the 
permutation $\sigma.$ We denote $\vert \sigma\vert$ the integer product of the cardinal of orbits of the permutation   $\sigma.$
One defines  $F_{(P,\pi,\sigma,\lambda_\pi)}:Im\Lambda_P\rightarrow Im\Lambda_P, \lambda_P\mapsto 
\sigma(\lambda_P)/ (\lambda_P\sigma(\lambda_\pi))$, 
the set $X_{P,\pi,\sigma,\lambda_\pi}=F_{(P,\pi,\sigma,\lambda_\pi)}(Im\Lambda_P)\cap Im \Lambda_{P_\sigma}$ is finite and we denote 
$\{ \lambda_\pi^\sigma\}\subset Im\Lambda_P $ a  set such that the restriction 
$F_{(P,\pi,\sigma,\lambda_\pi)}:\{ \lambda_\pi^\sigma\}\rightarrow X_{P,\pi,\sigma,\lambda_\pi}$ is a bijection.
In particular we have $\sigma(\lambda_\pi^\sigma)/\lambda_\pi^\sigma\sigma(\lambda_\pi)\in Im \Lambda_{P_{\sigma}}$ i.e is 
fixed by $\sigma.$
Note that the operator $M_{P,\tau\sigma}^{\tau(P)}(.,\lambda\lambda_\pi^\sigma/\mu_\sigma)$  and  $
[\tau\sigma(\lambda_\pi)]\circ M_{P,\tau\sigma}^{\tau(P)}(.,\lambda\lambda_\pi^\sigma)$ are vector space isomorphisms  from  
$L^2(M_P(F)N_P(\aa)\bc G(\aa)/J,\pi)$ to 
 $L^2(M_{\tau(P)}(F)N_{\tau(P)}(\aa)\bc G(\aa)/J,\tau(\pi))$ 
 because $\tau\sigma(\pi)=\tau(\pi)\otimes \tau\sigma(\lambda_\pi)^{-1}.$

\subsection{The simple spectral side}

If $S$ is a finite set of places of $F$, we will write $G_S$ for the Cartesian product $\times_{v\in S} G_v$, and $G^S$ for the restricted product $\t'_{v\notin S} G_v$. The same, if  $f\in H(\tg^D(\aa))$, we write $f_S$ for $\otimes_{v\in S} f_v$ (viewed as a function on $G_S$) and $f^S$ for $\otimes_{v\notin S} f_v$ (view as a function on $G^S$). If $Q$ is a standard parabolic subgroup of $G$ with Levi decomposition $Q=MN$, we adopt the same notation $Q_S$, $M_S$ etc.. Recall the definition of the constant term along the parabolic group $Q_S$ of a function $f_S$ like before: it is the function $f_S^Q$ defined on $M_S$ by the formula 
$$f_S^Q(l):=\delta_{Q_S}^{1/2}(l)\int_{N_S}\int_{K_S}f_S(k^{-1}lnk)dkdn$$
for every $l\in M_S$, where $\delta_{Q_S}$ is the modulus function of $Q_S$ (which plays no role here as we will show and use only that the integral vanishes under particular hypothesis). If $S=\{v\}$, i.e. is made of only one place, we replace index $S$ simply by index $v$.

The subgroup $J$ of $G(\aa)$ is not product. However, by choice of the generator $a$ of $J$, we have that $J$ is isomorphic to a subgroup of $G_{T_a}$ which we denote $J_{T_a}$, and we see $G(\aa)/J$ as the product $G_0\t G_V\t G^{T_a\cup V}$, where $G_0=G_{T_a}/J_{T_a}$ (recall we chose $T_a$ disjoint from $V$). We use the same notation for $G'$, so that $G'(\aa)/J=G_0\t G'_V\t G^{T_a\cup V}$. 

We show a simple form of the spectral side of the trace formula for functions in $H(\tg^D(\aa))$. 

Let $f\in H(\tg^D(\aa))$,  and
set $h(g):=\sum_{z\in J} f(zg)$ (for each $g$ the sum is finite as the support of $f$ is compact). We see also $h$ as a map from  $G(\aa)/J$ to $\cc$ locally constant with compact support. Moreover, and this is important in the sequel, $h$ is a tensor product, namely $h=h_0\otimes (\otimes_{v\notin T_a} h_v)$ where $h_0$ is a function on the quotient group $G_0$ and, for $v\notin T_a$, we have $h_v=f_v$. 
\begin{prop}\label{simplespectral}
We have:
$$Tr^T(h)=\sum_{\pi} \tr\pi(h)$$
where $\pi$ runs over the set of discrete series of $G(\aa).$ 
\end{prop}

{\bf Proof.} We want to prove that the terms $Tr^{T}_{(P,\pi)}(h)$ associated to proper parabolic subgroups
 $P(\aa)$ in the Lafforgue's Theorem \ref{lafforgue} 
vanish  for functions $h$ as in the proposition. 
This will be implied by the vanishing of 
 $m^T_{(P,\pi,\sigma,\lambda_\pi,\lambda_\pi^\sigma)}(\lambda_\sigma,h)=
 Tr(M_{(P,\pi,\sigma,\lambda_\pi,\lambda_\pi^\sigma)}^T(.,\lambda_\sigma, h))$ for all 
$(P,\pi,\sigma,\lambda_\pi,\lambda_\pi^\sigma)$ and  $\lambda_\sigma\in Im_{\Lambda_{P_{\sigma}}}.$

We say $(P,\pi)$ is {\bf regular} if $Fix(P,\pi)$ is reduced to one single element, the identity.

In order to simplify the argument we first explain the vanishing of this term when $(P,\pi)$ is regular. 
This implies that  we have $P_{\sigma}=P$ and $\{ \lambda_\pi^\sigma\}$ can be chosen to be the  singleton 
$\{ 1\}$;  we set $\lambda:=\lambda_\sigma.$
Therefore 
$Tr^{T}_{(P,\pi,\sigma=1,\lambda_\pi)}(h)$ is given by the formula of proposition \ref{lafforgue} 
 and the   the expression  giving $\sum_{\lambda_\pi}M_{(P,\pi,id,\lambda_\pi,1)}^T(.,\lambda,h)$ is 
 exactly the formula (11.4.10) of Laumon \cite{Lau}.
 
Let $M(\aa)$ be a proper Levi subgroup of $G(\aa)$, the proof is the same as  the series 
of results contained in 11.5 to  11.8 in \cite{Lau} which apply as soon as $\pi$ is regular,  based themselves on results
 of Arthur and particularly  splitting formula for
 $(G,M)$ families 
(see for example \cite{Ar1} and \cite{Ar2}).

Let $M$ be a standard Levi subgroup of $G$. Let $(n_1,n_2,...,n_k)$ the partition of $n$ associated to $M$. We say that $M$ transfers at the place $v\in V$ if $d_v|n_i$ for all $1\leq i\leq k$ (recall $V$ is the set of places where $D$ does not split).

\begin{lemme}\label{nontransfer}
If $M$ is proper, then there are at least two places in $V$ where $M$ does not transfer.
\end{lemme}
\ \\
{\bf Proof.} This comes from arithmetic consideration. For $v\in V$, we have $G'_v=GL_{r_v}(D_v)$ where $\dim_{F_v} D_v=d_v^2$ and 
 $r_vd_v=n$. 
  According to class field theory, we have that the Hasse invariant of $D$ at 
any place $v\in V$ is of 
the form $\frac{r_v x_v}{n}$,  with $x_v$ positive integer and $gcd(x_v,d_v)=1$. Moreover, $\frac{x_v}{d_v}$ is the Hasse invariant of $D_v$ and $\sum_v\frac{r_v x_v}{n}$is an integer, which we prefer to write as: {\it $n$ divides $\sum_{v\in V} r_v x_v$}. This is true in case $A$ is a simple central algebra over $F$. As here $A$ is, moreover, a division algebra, the least common multiple of $d_v, v\in V$ is $n$. 

Let $m$ be the greatest common divisor of the $n_i$. As $M$ is proper, $m<n$, and, as the least common multiple of $d_v, v\in V$ is $n$, there exists at least one place, $v_0\in V$, such that $d_{v_0}$ does not divide $m$. So $M$ does not transfer at $v_0$. But we also know that $n$ divides $\sum_{v\in V} r_v x_v$.

 Assume, for every $v\in V$, $v\neq v_0$, we had $d_v|m$. As $n=r_vd_v$, one has $n|mr_v$ for every $v\in V$, $v\neq v_0$. As 
$n|\sum_{v\in V} r_v x_v$, we have $n|mr_{v_0}x_{v_0}$. Then $d_{v_0}| m x_{v_0}$. But $gcd(d_{v_0},x_{v_0})=1$, so $d_{v_0}|m$ which is not possible.\qed

We will use this lemma to simplify the trace formula: one of the two places, so to say, for applying the splitting formula of Arthur (Laumon for characteristic $p$), the other one to kill almost all the remaining terms. 

 \begin{lemme}\label{vanishinglemma}
  Let $M$ be a proper standard Levi subgroup of $G$ and $v_0$ a place where $M$ does not transfer. Set $V_\infty:=V\bc \{v_0\}$. Let $f\in H(\tg^D(\aa))$,  then we have:
 
 (a) For every proper parabolic subgroup $Q$ of $G$ containing $M$, $f_{V_\infty}(k^{-1}xk)=0,$ for all $k\in K_{V_\infty}$, and all $x\in Q_{V_\infty}$. In particular, 
 $$f_{V_\infty}^{Q_{V_\infty}}=0.$$
 
 (b)  $Tr(I_{P^{V_\infty}}^{G^{V_\infty}}(\pi^{V_\infty}(\lambda)(f^{V_\infty})))=0$. 
 \end{lemme}
\ \\
{\bf Proof.} (a) Let $L$ be the Levi component of $Q$ containing $M$. According to Lemma \ref{nontransfer}, there are two places in $V$ where $L$ does not transfer, so at least one place $v_L$ in $V_\infty$. The support of $f_{v_L}$ contains solely 
 elements $g$ such that $P_g$ has irreducible factors of degree all divisible by $d_{v_L}$. Any element in $Q_{v_L}$ has characteristic polynomial which is a product of polynomials of degrees equal to the sizes of the blocs of $L$. So, no element in $Q_{v_L}$ may be conjugated to an element in the support of $f_{v_L}$. Now $f_{V_\infty}$ is a tensor product of functions, one of which is $f_{v_L}$ and the result follows.
 
 (b) The global trace is a product of local traces, and it is enough to prove that $Tr(I_{P_{v_0}}^{G_{v_0}}(\pi_{v_0})(\lambda)(f_{v_0})))=0$. By the same argument as before, because $M$ does not transfer at the place $v_0$, the support of $f_{v_0}$ does not meet any conjugated of $M_{v_0}$, and the constant term of $f_{v_0}$ along a parabolic subgroup having $M_{v_0}$ as Levi component vanishes. The Lemma 7.5.7 of Laumon \cite{Lau}(Part I page 189) shows then that the trace of the induced representation vanishes on $f_{v_0}$.\qed
 
\smallskip

 We now apply the series of results contained in  Laumon \cite{Lau}  which are based on the notion of $(G,M)$ family.
  The properties of $(G,M)$ families and of the weighted mean values  have been first introduced by Arthur and their 
  definition and 
  properties are  recalled in the review article of
 Arthur \cite{Ar}. 
 
 Let us recall that when $(c_P)_P$ is a $(G,M)$ family of holomorphic functions on $\Lambda_P,$ one can associate to it the function $c_M$ 
 (called in the sequel {\it the
 weighted  mean value} of $(c_P)$) defined by the proposition 11.5.7
 of Laumon \cite{Lau} i.e
 $$c_M(\mu)=\sum_{P\in {\mathcal P}(M)}\frac{c_P(\mu)}{\theta_P(\mu)}$$
 where $\theta_P(\mu)=\prod_{\alpha\in \Delta_P}(1-\check{\alpha}(\mu))$ and ${\mathcal P}(M)$ denotes the set of parabolic subgroups having $M$ as Levi component.
 The meromorphic function $c_M$ admits a holomorphic extension on the whole $\Lambda_P.$

 One can define a  restriction operator on $(G,M)$ families recalled in \cite{Ar}:
 if $Q$ is  a parabolic subgroup of $G$ containing $M$, we denote $M_Q$ its Levi subgroup, to each $(G,M)$  family $c$ one associates a
  $(M_Q,M)$ family denoted  $c^{Q}$.
  One has the splitting formula of  Arthur which enables to evaluate the weighted mean value of the product of two $(G,M)$ families 
  $c,c'$ as:
  $$
 (cc')_M=\sum_{L_1,L_2\in {\mathcal L}(M)} d_M^G(L_1,L_2) c_{L_1}^{P_{L_1}} {c'}_{L_2}^{P_{L_2}}.
 $$
 where ${\mathcal L}(M)$ is the set of Levi components of parabolic subgroups containing $M$, $P_L$
 is a certain parabolic subgroup of $G$ having $L$ as Levi component and $d_M^G(L_1,L_2)$ are
 complex coefficients which definition are recalled in \cite{Ar}.
 
 The properties of $(G,M)$ families and of the weighted mean values  are recalled in the review article of
 Arthur \cite{Ar}.

Laumon follows  three major steps:
 
1) One first expresses $m^T_{(P,\pi,\sigma=id,\lambda_\pi,1)}(\lambda,h)$ as the value at $\lambda_\pi^{-1}$ of the weighted  mean value of 
  a product of  two $(G,M)$ families as in example 11.5.9 of \cite{Lau}.

Indeed in the present notations (we denote $\mu=\mu_{\sigma},)$ we can define 
 $(G,M)$ families $c(\lambda_\pi;.), c'(.)$ as follows:

$$c_{\tau}(\lambda_\pi;\mu)= Tr_{L^2(M_P(F)N_P(\aa)\bc G(\aa)/J,\pi)} (([\tau(\lambda_\pi)]\circ M_{P,\tau}^{\tau(P)}(.,\lambda))^{-1}\circ M_{P,\tau}^{\tau(P)}{(.,\frac{\lambda}{\lambda_\pi\mu})
\circ h(.,\frac{\lambda}{\lambda_\pi\mu})}).$$ This  is a $(G,M)$ family of functions of $\mu$ indexed by 
$\tau\in 
\mathfrak{S}_{\vert P\vert}$ (the set of parabolic group ${\mathcal P}(M)$ which normally indexes 
the $(G,M)$ family is here indexed by  $\tau$ because  
${\mathcal P}(M)=\{\tau(P) ,\tau \in 
\mathfrak{S}_{\vert P\vert}\}$), and 
$c'{}_{\tau}(\mu)=
{\hat{{\bf 1}}}_{P,\tau}^T(\mu) \theta_\tau(\mu),$  with  $\theta_\tau=\theta_P$ (where $P$ corresponds to $\tau$)
and we have 
$$m^T_{(P,\pi,\sigma=id,\lambda_\pi,1)}(\lambda,h)=\lim_{\mu\rightarrow 1}(c(\lambda_\pi;.) c'(.))_{M}
 (\mu\lambda_{\pi}^{-1}),$$ 
this is exactly the formula Example (11.5.9) of \cite{Lau}.

Remark: $c'$ is a $(G,M)$ family as soon as $T\in {\mathfrak a}_{\emptyset,\mathbb{Z}}$ where ${\mathfrak a}_{\emptyset,\mathbb{Z}}
\subset {\mathfrak a}_{\emptyset}$ is the
root lattice generated by $(\alpha_i).$ We assume in the sequel,  as in \cite{Lau},  that $T$ satisfies this integrality condition.

2) Because $h$ satisfies the fact that its constant term $h^{Q}$ vanishes for every proper
parabolic subgroup containing $M(\aa)$, we have that the weighted mean value of the $(M_Q,M)$
family 
$c(\lambda_\pi,.)^Q$ is equal to $0$ (lemma 11.5.15 \cite{Lau}),
therefore using the splitting formula of Arthur we obtain that 
$m^T_{(P,\pi,id,\lambda_\pi,1)}(\lambda,h)=
c_M ( \lambda_\pi;\lambda_\pi^{-1}) c'_{P}(\lambda_\pi^{-1}).$ 
From  the analysis of Laumon (corollary 11.5.14 \cite{Lau}) if $c_M ( \lambda_\pi;\lambda_\pi^{-1})$ is different from zero 
then $\lambda_\pi$ is the restriction to $M$ of a character of $G$.
  Moreover by an explicit computation involving the exact 
 expression of $c'_{P}$ given by the formula of section 10.1 of \cite{NDT} or by the equivalent expression (Lemma 11.5.5 ii) of
 \cite{Lau} one obtains that when $\lambda_\pi$ is the restriction to $M$ of a character of $G$, $c'_{P}(\lambda_\pi^{-1})=0$ unless 
 $\lambda_\pi=1$ and in this case $c'_{P}(\lambda_\pi^{-1})=\vert \Gamma_P\vert$ where $\Gamma_P$ is the finite group appearing in 
 \cite{NDT}.$\vert \Gamma_P\vert$ can be computed and is equal  to $e=n/m.$

Therefore $m^T_{(P,\pi,id,\lambda_\pi,1)}(\lambda,h)$ is null unless $\lambda_\pi$ is trivial and in this case we have 

$$
m^T_{(P,\pi,id,\lambda_\pi=1,1)}(\lambda,h)=\vert \Gamma_P \vert
Tr({\mathcal R}_{M}(\pi,\lambda)\circ h(.,\lambda))
$$
where one defines the weighted mean value operator 
$${\mathcal R}_{M}(\pi,\lambda)=\lim_{\mu\rightarrow 1}
 \sum_{\tau\in  \mathfrak{S}_{\vert P\vert} }\frac{1}{\theta_\tau(\mu)}
(M_{P,\tau}^{\tau(P)}(.,\lambda))^{-1}
\circ M_{P,\tau}^{\tau(P)}(.,\lambda/\mu).$$

3) Using this proposition, one can  then express  $m_{P,\pi,id,\lambda_\pi=1,1}^T(\lambda,h)$,  as the value at $\mu=1$ of the  weighted
 mean value of the 
product of two 
$(G,M)$ families $c_{V_\infty}, c^{V_\infty}$ given by straightforward generalization of the  Lemma 11.6.6 of  \cite{Lau} (one has to replace $\infty$ in his
formulas by the finite set $V_\infty$). We assume that $h=h_{V_\infty}\otimes h^{V_\infty}$. The 
 vanishing of the constant term of $h_{V_\infty}$ for every proper parabolic subgroup implies that
 by the splitting formula we have the
factorization given by proposition 11.6.8 of \cite{Lau}:
$$m_{P,\pi,id,1}^T(\lambda,h)=\vert \Gamma_P\vert Tr({\mathcal R}_{M}(\pi_{V_\infty},\lambda)
\circ Ind_{P(F_{V_\infty})}^{G(F_{V_\infty})}(\pi_{V_\infty}(\lambda)({\check h}_{V_\infty})) 
Tr(I_{P^{V_\infty}}^{G^{V_\infty}}(\pi^{V_\infty}(\lambda)({\check h}^V_\infty))),$$
and
 ${\mathcal R}_{M}(\pi_{V_\infty},\lambda)$ is the generalization of the operator given by Laumon page 194 \cite{Lau} which reads in our case:
$${\mathcal R}_{M}(\pi_{V_\infty},\lambda)=
\lim_{\mu\rightarrow 1}
 \sum_{\tau\in  \mathfrak{S}_{\vert P\vert} }\frac{1}{\theta_\tau(\mu)}(\widehat{M}_{P,\tau}^{\tau(P)}(.,\pi_{V_\infty},\lambda))^{-1}
\circ \widehat{M}_{P,\tau}^{\tau(P)}(.,\pi_{V_\infty},\lambda/\mu),$$
where  $\widehat{M}_{P,\tau}^{\tau(P)}(.,\pi_{V_\infty},\lambda)=\bigotimes_{v\in V_{\infty}} 
\widehat{M}_{P,\tau}^{\tau(P)}(.,\pi_{v},\lambda)$ are tensor product  of the normalized local intertwining operator of Langlands-Shahidi, see theorem 11.6.4
of  \cite{Lau}.

After these steps which give an explicit expression for  $m^T_{(P,\pi,\sigma,\lambda_\pi)}(\lambda,h)$ in term of
 local components, it is sufficient to show that  $  Tr(I_{P^{V_\infty}}^{G^{V_\infty}}(\pi^{V_\infty}(\lambda)({\check h}^{V_\infty}))) $ vanishes, which holds as a consequence of Lemma (\ref{vanishinglemma} ).

We therefore have shown that each of the term $Tr^{T}_{(P,\pi,\sigma=id,\lambda_\pi)}(h)$ vanishes when $P$ is proper and $\pi$ is regular.

\medskip

When $\pi$ is not regular we can generalize the previous arguments as follows.

We fix $(\sigma,\lambda_\pi)\in Fix(P,\pi)$ and we fix a choice of $\{\lambda_{\pi}^\sigma\}.$ 

Step 1. amounts to show that  
$
m_{(P,\pi,\sigma,\lambda_\pi,\lambda_\pi^\sigma )}^T(\lambda_\sigma,h)$
is the evaluation at $\frac{\sigma(\lambda_\pi^\sigma)}{\lambda_\pi^\sigma \sigma(\lambda_\pi)}$ of the weighted mean value 
of the product of two $(G,M_\sigma)$ families.

We can define  $(G,M_\sigma)$ families of functions, (this is proven in \cite{NDT} 
proposition 10.8, lemma 11.9 and corollary 11.10),  $c(\lambda_\pi,\lambda_\pi^\sigma ;.), c'(.)$ on $ \Lambda_{P_{\sigma}}$ as:
\begin{eqnarray*}
&&c_{\tau}(\lambda_\pi,,\lambda_\pi^\sigma ;\mu_\sigma)=\\ 
&&Tr_{L^2(M_P(F)N_P(\aa)\bc G(\aa)/J,\pi)} 
([\tau\sigma(\lambda_\pi)]\circ (M_{P,\tau\sigma}^{\tau(P)}(.,\lambda_{\sigma}\lambda_{\pi}^\sigma))^{-1}\circ M_{P,\tau}^{\tau(P)}(.,\frac{\lambda_\sigma \sigma(\lambda_\pi^\sigma)}{\mu_{\sigma} \sigma(\lambda_\pi)})
\circ h(.,\frac{\lambda_\sigma \sigma(\lambda_\pi^\sigma)}{\mu_{\sigma} \sigma(\lambda_\pi)}))
\end{eqnarray*}
  and 
$c'{}_{\tau}(\mu_\sigma)=
{\hat{{\bf 1}}}_{P_\sigma,\tau}^T(\mu_\sigma) \theta_{\tau}(\mu_\sigma),$ 
where these  $(G,M_\sigma)$ families are  indexed by $\tau\in 
\mathfrak{S}_{\vert P_\sigma\vert}.$ 

We have 
$$m^T_{(P,\pi,\sigma,\lambda_\pi,\lambda_\pi^\sigma)}(\lambda_\sigma,h)=
\lim_{\mu_\sigma\rightarrow 1}(c(\lambda_\pi,\lambda_\pi^\sigma ;.) c'(.))_{M_\sigma}
 (\mu_\sigma \frac{\sigma(\lambda_\pi^\sigma)}{\lambda_\pi^\sigma \sigma(\lambda_\pi)}).$$ 
 
Step 2. can be modified as follows. 
 Because $h$ satisfies the fact that its constant term $h^{Q}$ vanishes for every proper
parabolic subgroup $Q$ containing $M$, we have that $h^{Q}$ vanishes for every proper parabolic containing $M_\sigma\supset M.$
As a result 
we have that the weighted mean value of the $(M_Q,M_{\sigma})$
family  $c(\lambda_\pi,\lambda_\pi^{\sigma}; .)^Q$ is equal to $0$.
As a result, using the splitting formula we obtain that 
$$m^T_{(P,\pi,\sigma,\lambda_\pi,\lambda_{\pi}^\sigma)}(\lambda_\sigma,h)=
c_{M_\sigma} (\frac{\sigma(\lambda_\pi^\sigma)}{\lambda_\pi^\sigma \sigma(\lambda_\pi)})
c'_{P_{\sigma}}(\frac{\sigma(\lambda_\pi^\sigma)}{\lambda_\pi^\sigma \sigma(\lambda_\pi)}).$$ 

We can do the same analysis as Laumon: the last expression is null unless 
$\frac{\sigma(\lambda_\pi^\sigma)}{\lambda_\pi^\sigma \sigma(\lambda_\pi)}$ is  the restriction to $M_\sigma$ of a character of $G.$
In this case,
 $c'_{P_{\sigma}}(\frac{\sigma(\lambda_\pi^\sigma)}{\lambda_\pi^\sigma \sigma(\lambda_\pi)})$ is null unless 
 $\frac{\sigma(\lambda_\pi^\sigma)}{\lambda_\pi^\sigma \sigma(\lambda_\pi)}$ is trivial.

Therefore we have the formula:
\begin{equation*}
m_{(P,\pi,\sigma,\lambda_\pi,\lambda_\pi^\sigma)}^T(\lambda_\sigma,h)=\vert \Gamma_{P_\sigma} \vert
Tr({\mathcal R}^\sigma_M(\pi,\lambda_\sigma\lambda_\pi^\sigma)\circ h(\cdot, \lambda_\sigma\lambda_\pi^\sigma)),
\end{equation*}

where $${\mathcal R}^\sigma_{M}(\pi,\lambda_\sigma\lambda_\pi^\sigma)=
\lim_{\mu_\sigma\rightarrow 1}\sum_{\tau\in  \mathfrak{S}_{\vert P_\sigma\vert} }\frac{1}{\theta_\tau(\mu_\sigma)}
([\tau\sigma(\lambda_\pi)]\circ  M_{P,\tau\sigma}^{\tau(P)}(.,\lambda_{\sigma}\lambda_{\pi}^\sigma))^{-1}\circ 
M_{P,\tau}^{\tau(P)}(.,\frac{\lambda_\sigma \lambda_{\pi}^\sigma}{\mu_{\sigma}})
$$

Step  3.  reduces to the fact  that the right hand side can be expressed as the weighted mean value at $\mu_\sigma=1$ of 
the  product of two 
$(G,M_{\sigma})$ families  defined by:

\begin{eqnarray*}
&&c_{\tau,V_{\infty}}(\mu_\sigma)=tr(R^\sigma_{\tau}(\mu_\sigma)\circ Ind_{P(F_{V_\infty})}^{G(F_{V_\infty})}(\pi_{V_{\infty}}(\lambda_\sigma\lambda_\pi^\sigma))({\check h}_{V_{\infty}}))\\
&&c_\tau^{V_{\infty}}(\mu_\sigma)=tr(S^\sigma_\tau(\mu_\sigma)\circ Ind_{P^{V_{\infty}}}^{G^{V_{\infty}}}(\pi^{V_{\infty}}(\lambda_\sigma\lambda_\pi^\sigma))({\check h}^{V_{\infty}}))
\end{eqnarray*}
where 

$$R^\sigma_{\tau}(\mu_\sigma)=\bigotimes_{v\in V_{\infty}} 
(\widehat{M}_{P,\tau\sigma}^{\tau(P)}(.,\pi_w,\lambda_{\sigma}\lambda_{\pi}^\sigma)
\circ[\tau\sigma(\lambda_\pi)])^{-1}\circ 
\widehat{M}_{P,\tau}^{\tau(P)}(.,\pi_w,\frac{\lambda_\sigma \lambda_{\pi}^\sigma}{\mu_{\sigma}}),$$
where  $\widehat{M}_{P,\tau}^{\tau(P)}(.,\pi_w,\lambda)$ are the normalized local intertwining 
operator of Langlands-Shahidi defined by 
$\widehat{M}_{P,\tau}^{\tau(P)}(.,\pi_w,\lambda)=a_\tau(\pi_w,\lambda){M}_{P,\tau}^{\tau(P)}(.,\pi_w,\lambda)$, 
${M}_{P,\tau}^{\tau(P)}(.,\pi_w,\lambda)$ is the local part at place $w$ of ${M}_{P,\tau}^{\tau(P)}(.,\pi,\lambda)$
and $a_\tau(\pi_w,.)$ are rational functions of the variable $\lambda\in \Lambda_P,$ whose properties are recalled in the Theorem 11.6.4 of 
 \cite{Lau}.
 \begin{eqnarray*}
&& S^\sigma_{\tau}(\mu_\sigma)=\prod_{v\in V_{\infty}} (a_\tau(\pi_w,\lambda_{\sigma}\lambda_{\pi}^\sigma)^{-1} 
a_\tau(\pi_w,\frac{\lambda_\sigma \lambda_{\pi}^\sigma}{\mu_{\sigma}}))\\
&&\times([\tau\sigma(\lambda_\pi)]\circ {M}_{P,\tau\sigma}^{\tau(P)}(.,\pi^{V_\infty},\lambda_{\sigma}\lambda_{\pi}^\sigma))^{-1}\circ 
M_{P,\tau}^{\tau(P)}(.,\pi^{V_\infty}, \frac{\lambda_\sigma \lambda_{\pi}^\sigma}{\mu_{\sigma}}).
\end{eqnarray*}

One has to show that $c_{V_{\infty}}, c^{V_{\infty}} $ are two $(G,M_\sigma)$ families. $c_{V_{\infty}}$ is easily shown to be a $(G,M_\sigma)$ family, the 
only non trivial point, as in the regular
case, is to show that $c^{V_{\infty}} $ is also a $(G,M_\sigma)$ family.
Proving this  goes along the same line as the proof of \cite{Lau} Lemma 11.6.6. 
We assume that $h=h_{V_\infty}\otimes h^{V_\infty}$. The 
 vanishing of the constant term of $h_{V_\infty}$ for every proper parabolic subgroup containing $M_\sigma$  implies that
  one obtains the exact analog of the factorization formula which reads:
 
 \begin{eqnarray*}
&& \frac{1}{ \vert \Gamma_{P_\sigma}\vert }m_{(P,\pi,\sigma,\lambda_\pi,\lambda_\pi^\sigma)}^T(\lambda_\sigma,h)=
 (c_{V_{\infty}}, c^{V_{\infty}})_{M_\sigma}(1)\\
 &&=(c_{V_{\infty}})_{M_\sigma}(1)(c^{V_{\infty}})_{P_\sigma}(1)\\
 &&=Tr({\mathcal R}^\sigma_{M}(\pi_{V_{\infty}},\lambda_\sigma\lambda_{\pi}^\sigma)\circ
 Ind_{P(F_{V_{\infty}})}^{G(F_{V_{\infty}})}(\pi_{V_{\infty}}(\lambda_\sigma\lambda_{\pi}^\sigma)(\check{h}_{V_{\infty}})) \\
&&\times Tr({\mathcal S}^\sigma(\pi^{V_\infty},\lambda_\sigma\lambda_{\pi}^\sigma)  I_{P^{V_\infty}}^{G^{V_\infty}}(\pi^{V_\infty}(\lambda_\sigma\lambda_{\pi}^\sigma)(\check{h}^{V_\infty}))),
\end{eqnarray*}

where ${\mathcal R}^\sigma_{M}(\pi_{V_\infty},\lambda)$ is the operator  which reads in our case:
$${\mathcal R}^\sigma_{M}(\pi_{V_\infty},\lambda)=
\lim_{\mu_\sigma\rightarrow 1}
 \sum_{\tau\in  \mathfrak{S}_{\vert P_\sigma\vert} }\frac{1}{\theta_\tau(\mu_\sigma)}
 \bigotimes_{v\in V_\infty}
 ({\widehat{M}}_{P,\tau}^{\tau(P)}(.,\pi_w,\lambda))^{-1}
\circ {\widehat{M} }_{P,\tau}^{\tau(P)}(.,\pi_w,\lambda/\mu_\sigma).$$
and 
\begin{eqnarray*}
&&{\mathcal S}^\sigma(\pi^{V_\infty},\lambda_\sigma\lambda_{\pi}^\sigma)=S^\sigma_{\tau=id}(1)\\
&&=([\sigma(\lambda_\pi)]\circ {M}_{P,\sigma}^{P}(.,\pi^{V_\infty},\lambda_{\sigma}\lambda_{\pi}^\sigma))^{-1}\circ 
M_{P,id}^{P}(.,\pi^{V_\infty}, \lambda_\sigma \lambda_{\pi}^\sigma)=
\\
&&=([\sigma(\lambda_\pi)]\circ {M}_{P,\sigma}^{P}(.,\pi^{V_\infty},\lambda_{\sigma}\lambda_{\pi}^\sigma))^{-1}.
\end{eqnarray*}
$[\sigma(\lambda_{\sigma}\lambda_{\pi}^\sigma)]\circ {M}_{P,\sigma}^{P}(.,\pi^{V_\infty},\lambda_{\sigma}\lambda_{\pi}^\sigma)
\circ[\lambda_{\sigma}\lambda_{\pi}^\sigma]^{-1}$ is an intertwining operator between the representation 
$I_{P^{V_\infty}}^{G^{V_\infty}}(\pi^{{V_\infty}}(\lambda_\sigma\lambda_{\pi}^\sigma))$ and the representation 
$I_{P^{V_\infty}}^{G^{V_\infty}}(\sigma(\pi^{{V_\infty}}(\lambda_\sigma\lambda_{\pi}^\sigma)).$  Because 
 $\sigma(\pi\otimes \lambda_\sigma\lambda_{\pi}^\sigma)=\pi\otimes \lambda_\sigma\lambda_{\pi}^\sigma,$ due to  
 $\sigma(\lambda_\pi^\sigma)=\lambda_\pi^\sigma \sigma(\lambda_\pi),$
 we therefore have that  $[\lambda_{\sigma}\lambda_{\pi}^\sigma]\circ [\sigma(\lambda_\pi)]\circ 
 {M}_{P,\sigma}^{P}(.,\pi^{V_\infty},\lambda_{\sigma}\lambda_{\pi}^\sigma)\circ [\lambda_{\sigma}\lambda_{\pi}^\sigma]^{-1}$ 
 is an intertwining operator of the representation
  $I_{P^{V_\infty}}^{G^{V_\infty}}(\pi^{{V_\infty}}(\lambda_\sigma\lambda_{\pi}^\sigma))$, which is irreducible because it is locally
 induced from irreducible unitary. As a result  $[\sigma(\lambda_\pi)]\circ {M}_{P,\sigma}^{P}(.,\pi^{V_\infty},\lambda_{\sigma}\lambda_{\pi}^\sigma)$ is a scalar operator and it is
 therefore sufficient to show that 
 
$Tr(I_{P^{V_\infty}}^{G^{V_\infty}}(\pi^{{V_\infty}}(\lambda_\sigma\lambda_{\pi}^\sigma)(\check{h}^{V_\infty})))=0.$
But this is implied by the previous lemma.

This ends the proof.

\qed

\subsection{The simple geometric side}

We show a simple form of the geometric side of the trace formula for functions in $H(\tg^D(\aa))$. 
Like in the previous subsection, we let $f\in H(\tg^D(\aa))$ and
set $h(g):=\sum_{z\in J} f(zg)$ which we see as a map from  $G(\aa)/J$ to $\cc$ locally constant with compact support. Here again, we have to play this game between $h$ and $f$ for the reason that $h$ is a function on $G(\aa)/J$ and it is not properly speaking a tensor product over places.

\begin{prop}\label{simplegeometric}
We have
$$Tr^{T}(h)=\sum_{o\in \tilde{O}_{G(F)}^D} vol(G_{\gamma_o}(F)\bc G_{\gamma_o}(\aa)/J) \sum_{z\in J}\Phi(f;z\gamma_o).$$
\end{prop}

{\bf Proof.}
Recall

\begin{equation*}
Tr^{T}(h)=\int_{G(F) \bc G(\aa)/J} dg \sum_{P\in { \mathcal P}_0^s} (-1)^{\vert P\vert -1}\sum_{\delta\in P(F) \bc G(F)} 
K_{h,P}(\delta g, \delta g) {\bf 1}_{P}^T(\delta g)
\end{equation*}
where 
\begin{equation*}
K_{h,P}(x,y)=\int_{N_P(\aa)}\sum_{\gamma\in M_P(F)} h(x^{-1}\gamma n y).
\end{equation*}

As $f \in H(\tg(\aa)^D)$, by Proposition \ref{poly} (b) we have that $K_{h,P}(x,x)$ is null for proper $P$. So
$$Tr^{T}(h)=
\int_{G(F)\bc G(\aa)/J}\sum_{\gamma\in \tg(F)^D} h(g^{-1}\gamma g)\, dg.$$ 

Moreover, using the claim (c) of the same proposition ($\tg(\aa)^D$ is stable under conjugation), we have:
$$Tr^{T}(h)=\int_{G(F)\bc G(\aa)/J}\sum_{\gamma\in \tg(F)^D} h(g^{-1}\gamma g)\, dg=
\int_{G(F)\bc G(\aa)/J}\sum_{\gamma\in \tg(F)^D}\sum_{z\in J} f(g^{-1}z\gamma g)\, dg.$$

We have

$$Tr^{T}(h)=\int_{G(F)\bc G(\aa)/J}\sum_{\gamma\in \tg(F)^D}\sum_{z\in J} f(g^{-1}z\gamma g)\, dg=$$
$$=\int_{G(F)\bc G(\aa)/J}\sum_{o\in \tilde{O}_{G(F)}^D}\sum_{\gamma\in o} \sum_{z\in J} f(g^{-1}z\gamma g)\, dg=$$
$$=\sum_{O\in \tilde{O}_{G(F)}^D}\int_{\gf\bc\ga/J}\sum_{\gamma\in o} \sum_{z\in J} f(g^{-1}z\gamma g)\, dg=$$
$$=\sum_{o\in \tilde{O}_{G(F)}^D}\int_{\gf\bc\ga/J}\sum_{t\in G_{\gamma_o}(F)\bc\gf} \sum_{z\in J} f(g^{-1}t^{-1}z\gamma_o tg)\, dg=$$
$$=\sum_{o\in \tilde{O}_{G(F)}^D} \int_{ G_{\gamma_o}(F)\bc\ga/J} \sum_{z\in J} f(g^{-1}z\gamma_o g)\, dg=$$
$$=\sum_{o\in \tilde{O}_{G(F)}^D} vol(G_{\gamma_o}(F)\bc G_{\gamma_o}(\aa)/J) \sum_{z\in J}\Phi(f;z\gamma_o).$$\qed

As in the proof of Deligne-Kazhdan simple trace formula, manipulations are allowed as for these elements $\gamma_o$ everything converges.\\

\subsection{Comparison with $G'(\aa)$}

Let $f\in H(\tg^D(\aa))$, $f'\in H(\tg'(\aa))$ such that $f \stackrel{\aa}{\leftrightarrow} f'$.  Let $DS$ be the set of irreducible subrepresentations of $R_G$ and $DS'$ the set of irreducible subrepresentations of $R_{G'}$.

\begin{prop}\label{equalityfdt} We have:
$$\sum_{\pi\in DS}\tr\pi(f)=\sum_{\pi'\in DS'}\tr\pi'(f').$$
\end{prop}

{\bf Proof.} Set $h(g):=\sum_{z\in J} f(zg)$, $h'(g):=\sum_{z\in J} f'(zg)$ and consider $h$ and $h'$ as  maps from  $G(\aa)/J$ and respectively $G'(\aa)/J$ to $\cc$, locally constant with compact support. It is enough to prove $\sum_{\pi\in DS}\tr\pi(h)=\sum_{\pi'\in DS'}\tr\pi'(h')$, as $\tr\pi(f)=\tr\pi(h)$ (by definition, $\pi(f)=\int_{G(\aa)}f\pi$ while $\pi(h)=\int_{G(\aa)/J} h\pi$ and the central character of $\pi$ is trivial on $J$).
Due to the hypothesis on $f$, the Propositions \ref{simplespectral} and \ref{simplegeometric} imply:

$$\sum_{\pi\in DS}\tr\pi(h)=\sum_{o\in \tilde{O}_{G(F)}^D}
 vol(G_{\gamma_o}(F)\bc G_{\gamma_o}(\aa)/J) \sum_{z\in J}\Phi(f;z\gamma_o).$$

The group $G'(F)\bc G'(\aa)/J$ is compact (\cite{Laf} III.6 Lemme 5 (ii)). So we have a similar formula:
$$\sum_{\pi'\in DS'}\tr\pi'(h')=\sum_{o\in \tilde{O}_{G'(F)}} vol(G'_{\gamma_o}(F)\bc G'_{\gamma_o}(\aa)/J) \sum_{z\in J}\Phi(f';z\gamma_o).$$
where $\{\gamma'_o\}$ is a system of representatives for $\tilde{O}_{G'(F)}$ such that $\gamma'_o\in O$ for all $o\in \tilde{O}_{G'(F)}$.

The Lemma \ref{ellipticglobal} establishes the unique characteristic polynomial preserving bijection between  $\tilde{O}_{G(F)}^D$ and $\tilde{O}_{G'(F)}$.
We have then equality term by term between the right hand member of these two equalities due to choices of measures and functions compatible with the local transfer.\qed

\subsection{End of the proof}

Now the proof goes the standard way, following ideas of Langlands. As this was usually applied in zero characteristic, we recall briefly the steps giving when needed the argument in non zero characteristic. 

Let $\pi\in DS$ be $D$-compatible. Let $U$ be the set of places $v$ of $F$ such that $v\notin T_a$, $G'_v$ splits (i.e. $v\notin V$) and $\pi_v$ is spherical. Let $U^c$ be the set of places of $F$ not in $U$, which is known to be a finite set. 
Let $DS_\pi$ be the subset of $DS$ made of representations $\tau$ such that $\tau_v\simeq \pi_v$ for all $v\in U$. Let  
$DS'_\pi$ be the subset of $DS'$ made of representations $\tau'$ such that $\tau'_v\simeq \pi_v$ for all $v\in U$. Then we have, for $f,f'$ as before:
\begin{equation}\label{flath}
\sum_{\tau\in DS_\pi}\tr\tau(f)=\sum_{\tau'\in DS'_\pi}\tr\tau'(f').
\end{equation}

This relation \ref{flath} is known to be a consequence of the Proposition \ref{equalityfdt} and the beautiful proof due to Langlands is now "standard" (it is detailed in the paper of Flath \cite{Fl2} for example).
The proof comes from the fact that an absolutely convergent sum of characters of non-isomorphic unitary spherical representations of $G^{U^c}$ is null if and only if the sum is void. 
This is based on the Satake isomorphism and abstract functional analysis and do not require
zero characteristic. 

By multiplicity one theorem (\cite{Sh}, \cite{P-S}), $DS_\pi=\{\pi\}$. Now we take $f_v=f'_v=1_{K_v}$ for all $v\in U$. Then $\tr\pi_v(f_v)=1$ for $v\in U$. So the relation \ref{flath} becomes:

$$\prod_{v\in U^c}\tr\pi_v(f_v)=\sum_{\tau'\in DS'_\pi}\prod_{v\in U^c}\tr\tau'_v(f'_v).$$

We know (\cite{BB}, \cite{Ba} Theorem 3.2) that the number of non isomorphic representations in $DS'_\pi$ is finite. As representations in $R_{G'}$ appear with finite multiplicity, the number of elements of $DS'_\pi$ is finite. As the number of representations involved in the equality is finite, we may switch from traces to characters: 

$$\prod_{v\in U^c}\chi_{\pi_v}(g_v)=\sum_{\tau'\in DS'_\pi}\prod_{v\in U^c}\chi_{\tau'_v}(g'_v)$$
whenever, for every $v\in U^c$, $g_v\lra g'_v$. By Theorem \ref{bhls}, and the hypothesis that $\pi$ is $D$-compatible, we may "transfer" characters from left to right. Writing $\lj_v$ for the Jacquet-Langlands local correspondence for unitary representations at the place $v$:

$$0=\epsilon\prod_{v\in U^c}\chi_{\lj_v(\pi_v)}(g'_v)+\sum_{\tau'\in DS'_\pi}\prod_{v\in U^c}\chi_{\tau'_v}(g'_v)$$
where $\epsilon$ is a sign (which appears from the local transfer Theorem \ref{bhls}). If we assume the linear independence of characters on groups $G'_v$ we have linear independence of characters for their product and we find there is just one $\tau'$ in $DS'_\pi$ and it verifies $\tau'=\lj_v(\pi_v)$ which is what we want. So let us give references for the linear independence in non zero characteristic. In \cite{JL} lemma 7.1 the linear independence of traces is proved, and the proof is independent of the characteristic. To pass from this result to the linear independence of characters it is enough to know the local integrability of characters. For groups like $G'_v$ (i.e. local inner forms of $GL_n$ in non zero characteristic) this is proved in \cite{Bacrelle} and \cite{Le2}.\\

On the other direction, to show the surjectivity, we start with $\pi'\in DS'$, let $U'$ be the set of places $v$ of $F$ such that $G'_v$ splits and $\pi'_v$ is spherical. We shortly come to a relation of the same type as \ref{flath} 
$$
\sum_{\tau\in DS_{\pi'}}\tr\tau(f)=\sum_{\tau'\in DS'_{\pi'}}\tr\tau'(f').
$$
where now $DS_{\pi'}$ and $DS'_{\pi'}$ are made of representations which have the same local component as $\pi'$ at places in $U'$. By the multiplicity one theorem, $DS_{\pi'}$ is void or contain a single representation. Again, the local independence of characters will show that, as $DS'_{\pi'}$ is not void, $DS_{\pi'}$ is not void neither and that the unique representation it contains is $D$-compatible. Then everything goes the same until the end of the proof.\qed




\section{Answer to two questions in \cite{LRS}}

\ \\
Here we answer two questions from \cite{LRS}. Only the second one is directly related to the main result of this paper. But the same question is related in \cite{LRS} also to the first, so we take the opportunity to answer it here too.\\

{\bf 1.} Let $F$ be a local field of non zero characteristic and set $G:=GL_n(F)$.
Let $\pi$ be a square integrable representation of $G$. Denote $z(\pi)$ the Zelevinsky dual of $\pi$. In \cite{LRS}, section 13.8, the authors ask the following question. Is there a function $f\in H(G)$ such that

(i) the orbital integrals of $f$ are null on regular semisimple elements which are not elliptic,

(ii) if $u$ is an irreducible unitary representation then $\tr u(f)\neq 0$ if and only if $u$ is isomorphic to $\pi$ or $z(\pi)$?

Such a function is known to exist if the characteristic of $F$ is zero. We give here the proof in non zero characteristic. It is known that the Paley-Wiener theorem (\cite{BDK}) allows one to construct a pseudocoefficient for $\pi$, i.e. $f\in H(G)$ such that 

- $\tr\tau(f)=0$ for all fully induced representation $\tau$ from any proper parabolic subgroup of $G$,

- $\tr\tau(f)=0$ for all tempered representation $\tau$ of $G$ such that $\tau$ is not isomorphic to $\pi$.

- $\tr\pi(f)=1$.

A detailed proof of the existence may be found in \cite{BaManuscripta} theorem 2.2. It is proved (op. cit. Lemme 2.4) that $f$ satisfies then the property (i). Let us explain why $f$ satisfies (ii). Let $u$ be an irreducible unitary representation of $G$ such that $\tr u(f)\neq 0$. Then, in the Grothendieck group of smooth representations of finite length of $G$, $u$ is a sum  $u=\sum_{i=1}^k s_i$ of standard representations $s_i$ all of which have the same cuspidal support as $u$. A standard representation is always tempered or fully induced from a proper parabolic subgroup. The reader will find definitions and proofs in \cite{DKV}, A.4.f. Now $\tr u(f)=\sum_{i=1}^k \tr s_i(f)$ so there is some $s_i$ which verifies $\tr s_i(f)\neq 0$. So one of the representations $s_i$ has to be $\pi$. So $u$ has the same cuspidal support as $\pi$, i.e. a Zelevinsky segment. 
According to the Tadi\'c classification of unitary representations of $GL_n$ (\cite{Ta}, any characteristic), $u$ is fully induced from a product of Speh representations twisted with some characters. As $\tr \tau(f)= 0$ for any fully induced representation $\tau$ from any proper parabolic subgroup of $G$, the product contains only one term and $u$ is a Speh representation. The cuspidal support of a Speh representation is easy to describe directly from its very definition (see \cite{Ta} for the definition), in particular it is easy to see that it has multiplicities unless $u$ is isomorphic to $\pi$ or $z(\pi)$. This finishes the proof.\\

Remark that the question of \cite{LRS} is asked in the Aubert (and non Zelevinsky) dual setting. But in \cite{Au} it is proved that the two duals differ by the sign $(-1)^k$ where $k$ is the number of cuspidal representations in the cuspidal support of $\pi$.\\

A formula for the orbital integrals of $f$ on the elliptic set is also conjectured in \cite{LRS} 13.8, which follows, in characteristic $p$, from Theorems 4.3 (ii) and 5.1 of \cite{BaManuscripta}.\\

{\bf 2.} The second question asked in \cite{LRS} is their Hypothesis 14.23. The authors explain in 14.24 that this Hypothesis would follow from the global Jacquet-Langlands correspondence. We confirm that the global Jacquet-Langlands correspondence, as stated and proved in our Theorem \ref{correspondence}, implies the Hypothesis in the way described in \cite{LRS} 14.23. As remarked by the authors, together with the results proved here in {\bf 1}, the Hypothesis implies then their Conjecture 14.21. Also the global Jacquet-Langlands correspondence simplifies their proof of the Theorem 14.12, as they do remark in the Remark 14.12. Indeed, let $D$ be a central global division algebra of degree $n^2$ over $F$ and $\pi'$ any discrete series of $D^\t$ which is Steinberg at one split place. Then $\pi'$ corresponds by Jacquet-Langlands to a discrete series $\pi$ of $GL_n$ which is Steinberg at the same place. Then $\pi$ is cuspidal because it has a local component which is square integrable. So $\pi$ is generic at every place. So $\pi'$ is generic at every place where $D$ splits.

\end{document}